\documentclass[11pt]{report}
\usepackage{geometry}                % See geometry.pdf to learn the layout options. There are lots.
\usepackage[parfill]{parskip}    % Activate to begin paragraphs with an empty line rather than an indent
\usepackage{graphicx}
\usepackage{amssymb}
\usepackage{amsmath}
\usepackage{epstopdf}
\usepackage{verbatim}
\newcommand{\R}{\mathbb{R}}
\newcommand{\E}{\mathbb{E}}	
\renewcommand{\P}{\mathbb{P}}

\newcommand{\Wd}{\mathfrak{W}_d}	% Williamson d-transform
\newcommand{\MDA}{\mathrm{MDA}}
\newcommand{\RV}{\mathrm{RV}}

\newcommand{\0}{\mathbf{0}}
\newcommand{\bx}{\boldsymbol{x}}
\newcommand{\bu}{\boldsymbol{u}}
\newcommand{\bX}{\boldsymbol{X}}	
	
\newcommand{\bU}{\boldsymbol{U}}
\newcommand{\bW}{\boldsymbol{W}}

\newcommand{\eqdis}{\stackrel{\mbox{\tiny d}}{=}}
\newcommand{\I}{1\!\!1}

\newtheorem{theorem}{Theorem}[chapter]
\newtheorem{lemma}[theorem]{Lemma}
\newtheorem{definition}[theorem]{Definition}
\newtheorem{proposition}[theorem]{Proposition}
\newtheorem{corollary}[theorem]{Corollary}
\newtheorem{example}[theorem]{Example}

\def\proofend{{\hfill $\square$}\bigskip}

\title{\textsc{Tail Properties of\\ Multivariate Archimedean Copulas}}
\author{Martin Larsson \\ \\ \\ \\ \\ \\ \\ \\ \\ \textsc{Master of Advanced Studies in Finance} \\ \\ \textsc{ETH Zurich and University of Zurich} \\ \\ Supervised by Dr.~Johanna Ne\v{s}lehov\'a \\ \\ \\}
\date{August 2008}                                           % Activate to display a given date or no date

\begin{document}
\maketitle

\newpage

\chapter*{Acknowledgment}

I would like to thank Dr.~Johanna Ne\v{s}lehov\'a for introducing me to the theory of Archimedean copulas, and for indicating the direction of investigation taken in this thesis. I am very grateful that she agreed to be my supervisor. Her comments, remarks and ideas have been crucial for the final result. Thanks are also due to Prof.~Paul Embrechts for his constant support and encouragement.

\newpage

\tableofcontents

\newpage

\chapter{Introduction and preliminaries}

Most problems in engineering and applied mathematics are motivated by some real world situation where one wants to make quantitative statements about the behavior and evolution of some system. Examples are mechanical engines, chemical processes, weather systems, or financial markets. These systems are rarely completely deterministic (in the sense that they can be described exactly by a deterministic model), but exhibit some randomness in their behavior. Sometimes these effects are negligible, but often it is necessary to take them into account. In some applications the influence of the random effects dominate over the deterministic structure, and a model that accurately captures the stochastic properties of the system is indispensable. If in such a situation the randomness is induced by different sources, or if several quantities are influenced differently by one source, a central issue is how to describe the dependence between these quantities; that the system behaves randomly does not mean that all components are independent.

Dependence between random variables can be very complicated, and not only in high dimensions. For example, an important but difficult issue is to find a relevant measure of the \emph{strength} of dependence. It is easy to find examples where a very common such measure, the correlation coefficient, does not give any information at all\footnote{Let $X$ and $Y$ be independent random variables, suppose $\E[Y]=0$, and set $Z:=XY$. Then $\mathrm{Cov}[X,Z]=\E[XZ]-\E[X]\E[Z]=\E[X^2]\E[Y]-\E^2[X]\E[Y]$ by independence. The covariance is thus zero, but $X$ and $Z$ are clearly far from independent. See Embrechts et al.~\cite{Embrechts/McNeil/Straumann:2002} for more on the intricacies of the correlation coefficient.}.

Different approaches have been taken to study the concept of dependence, and a variety of results exist in probability theory. In this thesis we shall focus on so-called \emph{copulas}, in particular a certain class of copulas called \emph{Archimedean copulas}. A copula is an object that allows to isolate the dependence structure from the marginal behavior of a random vector, which probably is the main reason why the concept is theoretically useful. It also turns out that arbitrary marginal distributions can be combined to a multivariate joint distribution using any specified copula. This makes copulas a useful tool in practical applications, where it can be relatively easy to model the marginal distributions, but much harder to get an idea of the dependence. Using copulas, these two problems can be separated. Before continuing, we take the opportunity to mention that the copula approach to multivariate modeling is not uncontroversial. See Mikosch \cite{Mikosch:2007} for a discussion of the drawbacks.

What is then a copula? Formally it is defined as follows:

\begin{definition} \label{def:copula}
A \emph{$d$-dimensional copula}, generically denoted $C$, is a $d$-dimensional distribution function with standard uniform marginal distributions\footnote{Note that the condition on the marginals implies that the support of any $d$-dimensional copula is contained in the unit cube $[0,1]^d$.}.
\end{definition}

It is not clear from the definition why copulas are useful for studying dependence. The following basic result due to Sklar \cite{Sklar:1959} links the definition of copulas to their role in dependence modeling. We omit the proof, which is fairly simple in the case when the marginal distribution functions are continuous, and more involved when they are not.

\begin{theorem} \label{thm:Sklar_df}
Let $H$ be a $d$-dimensional distribution function and $F_1,\ldots,F_d$ its one-dimensional marginals. Then there exists a copula $C$ such that
$$
H(x_1,\ldots,x_d) = C(F_1(x_1),\ldots,F_d(x_d)).
$$
$C$ is uniquely determined on $\mathrm{Ran}\ F_1\times\ldots\times\mathrm{Ran}\ F_d$, where $\mathrm{Ran}\ F_i$ is the range of $F_i$. Conversely, given a copula $C$ and univariate distribution functions $F_1,\ldots,F_d$, the above relationship always defines a $d$-dimensional distribution function $H$.
\end{theorem}

Sklar's theorem clarifies the role of the copula $C$ as carrier of the dependence structure, fully separated from the marginals. A result similar to Theorem \ref{thm:Sklar_df} holds for survival functions. Although completely analogous we state it explicitly, since in our context survival functions will occur frequently.

\begin{theorem} \label{thm:Sklar}
Let $\bar{H}$ be a $d$-dimensional survival function and $\bar{F}_1,\ldots,\bar{F}_d$ the corresponding one-dimensional marginal survival functions. Then there exists a copula $C$ such that
$$
\bar{H}(x_1,\ldots,x_d) = C(\bar{F}_1(x_1),\ldots,\bar{F}_d(x_d)).
$$
The copula is uniquely determined on $\mathrm{Ran}\ \bar{F}_1\times\ldots\times\mathrm{Ran}\ \bar{F}_d$. Conversely, given a copula $C$ and univariate survival functions $\bar{F}_1,\ldots,\bar{F}_d$, the above relationship always defines a $d$-dimensional survival function $\bar{H}$.
\end{theorem}

The object $C$ in Theorem \ref{thm:Sklar} is called the \emph{survival copula} of $H$ (or $\bar H$). It is easy to show that if $C_1$ is the (ordinary) copula of $H$, and $C_2$ is its survival copula, then
$$
(U_1,\dots U_d)\sim C_1 \qquad \mathrm{if\ and\ only\ if} \qquad (1-U_1,\dots,1-U_d)\sim C_2,
$$
where as usual the symbol $\sim$ stands for ``is distributed according to''. It is important to be clear about the distinction between a copula and the corresponding survival copula. This sometimes becomes confusing, especially when the tails of a distribution are related to the tails of the copula (or survival copula). For example, the dependence in the upper tail of a distribution $H$ is captured by the upper tail of its (ordinary) copula, but by the lower tail of its survival copula.

We note that Sklar's Theorem implies that given a joint distribution function $H$ and its marginals $F_1,\dots,F_d$, the copula can be expressed as
\begin{equation} \label{eqn:intr_1}
C(u_1,\dots,u_d) = H(F^{\leftarrow}(u_1),\dots,F^{\leftarrow}(u_d)).
\end{equation}
Here $F_i^{\leftarrow}$ is the \emph{generalized inverse} of $F_i$, defined as
$$
F_i^{\leftarrow}(u) :=  \inf\{x \in \R | F_i(x)\geq u\}, \qquad u\in [0,1].
$$
Not all distribution functions $F$ are invertible, but since they are non-decreasing the generalized inverse is always well-defined. The generalized inverse of the survival function $\bar F_i$ is defined as $\bar F_i^{\leftarrow}(u):=F_i^{\leftarrow}(1-u)$ for $u\in [0,1]$. Similarly as in (\ref{eqn:intr_1}), the survival copula can be extracted from $\bar H$ and $\bar F_1,\dots,\bar F_d$.

As already mentioned, this thesis is devoted to the class of Archimedean copulas. A lot of research has been done in this field, both concerning the general theory of copulas and Archimedean copulas in particular. Nelsen \cite{Nelsen:2006} contains a rigorous mathematical introduction, and McNeil, Frey and Embrechts \cite{McNeil/Frey/Embrechts:2005} gives a good account of how to apply copulas in (financial) risk management. In this thesis we will exploit some recent advances in the understanding of the structure of such copulas. We give an introduction to Archimedean copulas and derive the necessary results and tools. This is carried out in Chapter \ref{section:Arch_copulas}. In particular our aim is to emphasize how these copulas arise as the dependence structures of certain mixture models.

The subsequent analysis will focus on the extremal behavior of Archimedean copulas. More precisely, we are interested in questions like how to describe the dependence between the largest observations (in a suitable sense) of a multivariate sample, or the dependence between observations exceeding some high threshold. Another issue is to what extent large outcomes occur simultaneously. All these concepts will be made precise as they are investigated. Some standard tools for dealing with such issues are introduced in Chapter \ref{section:RV_EVT}, where we give a brief account of extreme value theory (EVT) and the theory of regularly varying functions.

Chapter \ref{section:Wd_tail}, which contains the main theoretical contribution of this thesis, concerns the asymptotic behavior of a certain integral transform that proves to be of central importance in the theory of Archimedean copulas. The results in this section can, and should, also be interpreted as properties of the tails of certain Beta mixtures.

The last part of the thesis, Chapters \ref{section:taildep}--\ref{section:threshold}, have the character of exercises in extreme value theory. It is shown how the results and concepts developed in Chapters \ref{section:Arch_copulas}--\ref{section:Wd_tail} make it straightforward to determine extremal properties of Archimedean copulas, such as coefficients of tail dependence, extreme value limits and limiting threshold copulas (notions that will be given precise definitions later on).

It is worth noting that we only are concerned with distributional properties of the random elements under considerations. We therefore only remark in passing that one always can find a probability space $(\Omega,\mathcal F, \P)$ where all random elements live. This is a consequence of Skorokhod's Representation Theorem. That said, there will be no further mention of probability spaces.

Now some remarks on notation and terminology. Capital letters ($X$, $Y$, $U$, $\dots$) denote random variables. Vectors are written in boldface---so, for instance, $\bX=(X_1,\ldots,X_d)$ is a random vector taking values in $\R^d$, and $\bx=(x_1,\ldots,x_d)$ is a $d$-dimensional deterministic vector. The letter $d$ will denote dimension, and is always an integer greater than or equal to two. Equalities, inequalities and arithmetic operations are to be interpreted componentwise. For example, $\boldsymbol x>\boldsymbol y$ means $x_i>y_i$ for all $i$.

With a slight abuse of language we will say that $\bar H$ is a \emph{survival function on the set $A\subseteq\R^d$}, when we mean that the corresponding random vector takes its values in $A$ with probability one. In particular, saying that $\bar H$ is a survival function on $\R^d_+:=\{\boldsymbol x \in \R^d \mid \boldsymbol x \geq \boldsymbol 0\}$ means that a random vector $\bX \sim \bar H$ satisfies $\bX\in\R^d_+$ a.s. Saying that $\bar H$ is a survival function on the \emph{interior} of $\R^d_+$ means the same thing, except that we now exclude probability mass on the boundary of $\R^d_+$. Note that such statements about $\bar H$ do \emph{not} mean that the function $\bx \mapsto \bar H(\bx)$ has its support confined to $\R^d_+$ or its interior. They do, however, imply that $\bar H$ is completely determined by its values in this set.

Several results, particularly in the last chapters, will concern \emph{weak convergence} of random variables, denoted $X_n\Rightarrow X$ as $n\to \infty$. This means that the corresponding distribution functions $F_n$ of $X_n$ converge weakly to the distribution function $F$ of $X$, i.e.~that $\lim_{n\to\infty} F_n(x) = F(x)$ for all $x$ where the limit $F$ is continuous. Weak convergence of random vectors, $\boldsymbol X_n \Rightarrow \boldsymbol X$, is defined analogously.

Finally, since distributions of mixtures, i.e.~products of two random variables, will be a central theme, conditioning arguments are crucial. This can be quite a delicate issue (see e.g.~Proschan and Presnel \cite{Proschan/Presnell:1998} for good illustrations of some of the pitfalls), which we nevertheless can circumvent using the following result:

\begin{lemma} \label{lemma:condexpindep}
Suppose that the $\R^n$-valued random variable $\boldsymbol X$ and the $\R^m$-valued random variable $\boldsymbol Y$ are independent, and let $h:\R^n\times\R^m\to\R$ be a measurable function such that $f(\boldsymbol x):=\E[h(\boldsymbol x,\boldsymbol Y)]$ is well-defined for all $\boldsymbol x\in \R^n$. Then
$$
\E[h(\boldsymbol X,\boldsymbol Y) | \boldsymbol X] = f(\boldsymbol X).
$$
In particular, if $X$ is scalar with distribution function $F$, then
$$
\P[X\boldsymbol Y>\boldsymbol z]  = \int_{\R} \P[\boldsymbol Y > t^{-1}\boldsymbol z]dF(t).
$$
\end{lemma}

\begin{proof}
The first part is Corollary 4.38 in Breiman \cite{Breiman:1968}. For the last part, let $h(x,\boldsymbol y):=\I_{\{\boldsymbol y>x^{-1}\boldsymbol z\}}$ and note that $\P[X\boldsymbol Y>\boldsymbol z]=\E[\P[\boldsymbol Y>X^{-1}\boldsymbol z | X]] = \E[\E[h(X,\boldsymbol Y)|X]]$. Since $f(x)=\P[\boldsymbol Y>x^{-1}\boldsymbol z]$, the result follows from the first part of the lemma.
\proofend\end{proof}

\chapter{Archimedean copulas} \label{section:Arch_copulas}

Archimedean copulas originally emerged as a family of copulas characterized by certain algebraic properties. This chapter is devoted to investigating their structure and characteristics. Section \ref{section:Arch_def_ex} defines Archimedean copulas and the so-called \emph{Archimedean generators}. Some important examples are given. We then continue with some theoretical background on \emph{$d$-monotone functions} (Section \ref{section:dmon}) and \emph{$\ell_1$-norm symmetric distributions} (Section \ref{section:l1nsd}). This theory is necessary to prove (and understand) the main result, Theorem \ref{thm:main}, which is given in Section \ref{section:Arch_char}. This result shows that Archimedean copulas can be \emph{characterized} as the dependence structure of the very well-behaved family of $\ell_1$-norm symmetric distributions. It also characterizes the relevant class of Archimedean generators. Further we discuss interpretations of the Archimedean dependence structure. Section \ref{section:Arch_compres} concludes with a result that clarifies the correspondence between Archimedean copulas, their generators and $\ell_1$-norm symmetric distributions. The theory presented in this chapter (except Section \ref{section:Arch_compres}) is based on McNeil and Ne\v{s}lehov\'a \cite{McNeil/Neslehova:2008}.

\section{Definition and examples} \label{section:Arch_def_ex}

\begin{definition} \label{def:Arch} (Archimedean copula and Archimedean generator)
A $d$-dimensional copula $C$ is called \emph{Archimedean} if it has the representation
$$
C(\boldsymbol u) = \psi(\psi^{-1}(u_1)+\ldots+\psi^{-1}(u_d)), \qquad \boldsymbol u \in \R^d_+,
$$
where $\psi:\R_+\to\R_+$ is an \emph{Archimedean generator}, i.e.~continuous, strictly decreasing on $\{\psi>0\}$ and satisfying $\psi(0)=1$ and $\lim_{x\to\infty}\psi(x)=0$. Here $\psi^{-1}$ is defined as $\psi^{-1}(u):=\inf\{x>0 | \psi(x)=u\}$.

If $\psi^{-1}(0)=\infty$ the generator is said to be \emph{strict}, otherwise \emph{non-strict}.
\end{definition}

The function $\psi^{-1}$ is sometimes called the \emph{pseudo-inverse} of $\psi$. Due to monotonicity and continuity of $\psi$, it coincides with the usual inverse of the restriction of $\psi$ to $\{\psi>0\}$. If $\psi$ is strict, this set is the whole of $\R_+$. Non-strict generators are not invertible, but $\psi^{-1}$ always is a right inverse of $\psi$, i.e.~$\psi\circ\psi^{-1}(u)=u$ for all $u\in[0,1]$.

Many authors define Archimedean copulas in terms of $\phi=\psi^{-1}$, and refer to this function as the generator. Although at first sight it is just a matter of convention, our definition turns out to be much more natural for studying the structure of these copulas. In this thesis, $\psi$ is thus the generator, and $\phi=\psi^{-1}$ is referred to as the \emph{inverse generator}.

An important observation is that $\psi$ defines a survival function on $\R_+$. This will be used later on.

\begin{example} \label{ex:arch_copulas}
The copulas $\Pi(\boldsymbol u) := u_1u_2\cdots u_d$ (the \emph{independence copula}) and $W(u_1,u_2):=(u_1+u_2-1)_+$ (the \emph{countermonotonicity copula}) are both Archimedean, with generators $\psi_\Pi(x) = e^{-x}$ and $\psi_W(x)=(1-x)_+$, respectively.

The copulas $C^{\mathrm{Cl}}_{\theta}(\boldsymbol u)$ generated by $\psi^{\mathrm{Cl}}_\theta(x)=(1+\theta x)^{-1/\theta}$, $\theta\geq -1$ (the \emph{Clayton copula}), and $C^{\mathrm{Gu}}_\theta(\boldsymbol u)$ generated by $\psi^{\mathrm{Gu}}_\theta(x)=\exp\{-x^{1/\theta}\}$, $\theta\geq 1$ (the \emph{Gumbel copula}), are two examples of parametric families of Archimedean copulas. They are among the most commonly encountered, and turn out to be particularly important in our context.
\end{example}

It is easy to check that the functions $\psi_\Pi$, $\psi_W$, $\psi^{\mathrm{Cl}}_\theta$ and $\psi^{\mathrm{Gu}}_\theta$ are all Archimedean generators. As for $W$ (in the bivariate case) and $\Pi$ (in any dimension), we even know that these functions are copulas\footnote{The marginals are clearly standard uniform. $\Pi$ is a product of univariate distribution functions and hence itself a distribution function (with independent marginals). For $W$, observe that $W(u_1,u_2)=\P[U\leq u_1, 1-U\leq u_2]$ for $U$ standard uniform.}. But we have not checked that the parametric families of claimed copulas are indeed copulas in the sense of Definition \ref{def:copula}. To do this we need conditions for when an Archimedean generator $\psi$ actually generates an Archimedean copula in a given dimension. In Section \ref{section:Arch_char} we will characterize these generators and derive a very useful stochastic representation for Archimedean copulas. First, however, we need some background.

\section{$d$-monotonicity and the Williamson $d$-transform} \label{section:dmon}

In the paper \cite{Williamson:1956}, R.~E.~Williamson studies the properties of real-valued functions having the following regularity property:

\begin{definition} \label{def:dmon}  ($d$-monotonicity)
A function $f:[0,\infty)\to \R$ is called \emph{$d$-monotone} for $d\geq 2$ if it is $d-2$ times differentiable in $(0,\infty)$ and the derivatives satisfy
$$
(-1)^kf^{(k)}(x)\geq0,\quad k=0,1,\ldots,d-2,
$$
and $(-1)^{d-2} f^{(d-2)}$ is non-increasing and convex in $(0,\infty)$. It is said to be \emph{completely monotone} if the above holds for all $d\geq 2$.
\end{definition}

Note that although $f^{(d-2)}$ not necessarily has a derivative, it is convex and hence absolutely continuous on $(0,\infty)$. For ease of notation, we will write $f^{(d-1)}$ for the \emph{right-continuous version} of the density of $f^{(d-2)}$.

Williamson \cite{Williamson:1956} gives a representation theorem for $d$-monotone functions in terms of a certain integral transform, whose basic properties he also derives. These results were cast into the form presented here by McNeil and Ne\v{s}lehov\'a \cite{McNeil/Neslehova:2008}. 

\begin{definition} \label{def:Wd} (Williamson $d$-transform)
Let $X$ be a non-negative random variable with distribution function $F$ and $d\geq2$ an integer. The \emph{Williamson $d$-transform} of $X$ (or $F$) is a real function on $[0,\infty)$ given by
\[ \Wd F(x) = \int_x^\infty\left(1-\frac{x}{t}\right)^{d-1}dF(t) = \left\{ \begin{array}{ll}
  \E\left[ \left(1-\frac{x}{X}\right)_+^{d-1}\right] & \mathrm{if \ }x>0,   \\
  1-F(0)&  \mathrm{if\ }x=0.
\end{array}\right.. \]
\end{definition}

\begin{proposition} \label{prop:Wd_prop} (Properties of the Williamson $d$-transform)
Let $F$ be a distribution function on $\R_+$. The Williamson $d$-transform $f=\Wd F$ satisfies the following properties:
\begin{itemize}
\item[(i)] (Right continuity) $f(x)\to f(0)=1-F(0)$ as $x\downarrow 0$.
\item[(ii)] (Order of decay) $x^kf^{(k)}(x)\to 0$ as $x\to\infty$ for $k=0,1,\dots,d-1$.
\item[(iii)] (Invertibility) The Williamson $d$-transform is invertible, and
$$
F(x) = \Wd^{-1} f(x) = 1-\sum_{k=0}^{d-1}\frac{(-1)^k}{k!}x^kf^{(k)}(x), \qquad x\geq 0.
$$
\end{itemize}
\end{proposition}

\begin{proof}
Part $(i)$ follows from the Dominated Convergence Theorem upon noting that
$$
\left(1-\frac{x}{X}\right)^{d-1}_+ = \left(1-\frac{x}{X}\right)^{d-1}\I_{\{X>x\}} \to \I_{\{X>0\}} \mathrm{\ a.s.},\qquad x\downarrow 0.
$$
Part $(ii)$ is shown by Williamson \cite{Williamson:1956}, Lemma 1 $(ii)$, and part $(iii)$ is Proposition 5 $(ii)$ in McNeil and Ne\v{s}lehov\'a \cite{McNeil/Neslehova:2008}.
\proofend\end{proof}

The following result is the main reason for studying the Williamson $d$-transform in connection with Archimedean copulas. Its relevance will be apparent later on.

\begin{theorem} \label{thm:dmon_repr} (Representation theorem)
Let $f$ be a real-valued function on $\R_+$. Then $f$ is $d$-monotone with $f(0)=p \in [0,1]$ and $\lim_{x\to\infty}f(x)=0$ if and only if $f = \Wd F$ for some distribution function $F$ on $\R_+$.
\end{theorem}

\begin{proof}
The proof can be found in McNeil and Ne\v{s}lehov\'a \cite{McNeil/Neslehova:2008}, Proposition 5 $(i)$. It is essentially a reformulation of Williamson's original proof in terms of distribution functions.
\proofend\end{proof}

In Williamson's paper, the $d$-transform is viewed in the context of real analysis as a mapping on a set of non-decreasing functions on the positive half-line. With the formulation here we view it in a probabilistic sense as a mapping on the set of distributions supported on the positive half-line. Using Lemma \ref{lemma:condexpindep} (for the second equality below), this probabilistic interpretation can be taken further by noting that, for $G(x)=(1-x)^{d-1}_+$,
$$
\mathfrak W_d F(x) = \int_0^\infty G(x/t) dF(t) = \E[\P[X>x/Y|Y]] = \P[XY>x],
$$
where $X$ and $Y$ are independent, $X\sim F$, and $Y\sim G$ and hence follows a $\mathrm{Beta}(1,d-1)$ distribution. The Williamson $d$-transform thus describes the relationship between non-negative random variables $X$ and certain Beta mixtures $XY$.

We conclude this section with a result that is \emph{the} key to analyzing the structure of Archimedean copulas.

\begin{proposition} \label{prop:dmon_sf}
Let $f$ be a real-valued function on $[0,\infty)$ and $\bar H$ given by
$$
\bar H(\boldsymbol x) = f(x_1 + \dots + x_d), \qquad \boldsymbol x \in \R^d_+.
$$
Then $\bar H$ is a survival function on the interior of $\R^d_+$ if and only if $f$ is $d$-monotone and satisfies $f(0)=1$ and $\lim_{x\to\infty}f(x)=0$.
\end{proposition}

\begin{proof}
See Proposition 2 in McNeil and Ne\v{s}lehov\'a \cite{McNeil/Neslehova:2008}. The proof consists in showing that $d$-monotonicity is a necessary and sufficient condition for the (possibly signed) measure induced by $\bar H$ to put non-negative mass on so-called \emph{rectangle sets}, i.e.~sets of the form $(a_1,b_1]\times\ldots\times(a_d,b_d]\subset \R^d_+$. Together with conditions on $\bar H$ and $f$ at the origin and infinity, this is equivalent to $\bar H$ inducing a probability measure. Note also that in our terminology, $\bar H$ being a survival function on the interior of $\R^d_+$ excludes it from assigning mass to the boundary of this set.
\proofend \end{proof}

\section{$\ell_1$-norm symmetric distributions} \label{section:l1nsd}

Although the relatively simple algebraic form (see Definition \ref{def:Arch}) is one of the reasons why Archimedean copulas are popular in some situations, this form does not help to interpret the corresponding dependence structure. Moreover, in many cases computations still become involved and seldom lead to any deeper understanding of \emph{why} results hold true. This is in particular true in high dimensions, where calculations become cumbersome. Instead of working directly with the copulas, we would therefore prefer, if possible, to work with some suitable class of distributions having a simple structure. The family of  \emph{$\ell_1$-norm symmetric distributions} introduced in this section turns out to be an appropriate class. In what follows, $\mathcal S^d$ will denote the \emph{$d$-dimensional unit simplex}:
$$
\mathcal S^d := \{ \bx \in \R^d_+ \ | \ \|\bx\|_1=1 \}, \qquad \mathrm{where} \quad  \|\bx\|_1=\sum_{i=1}^d |x_i|.
$$

\begin{definition}
A random vector $\boldsymbol X$ on $\R^d_+$ is called \emph{$\ell_1$-norm symmetric} if there exists a non-negative random variable $R$ and a random vector $\boldsymbol S$, independent of $R$ and uniformly distributed on the $d$-dimensional simplex $\mathcal S^d$, such that
$$
\boldsymbol X \eqdis R \boldsymbol S.
$$
The random variable $R$ is referred to as the \emph{radial part} of $\boldsymbol X$ and its distribution as the \emph{radial distribution}.
\end{definition}

The stochastic representation in terms of a certain simplex mixture makes $\ell_1$-norm symmetric distributions very tractable objects. In the following we restrict attention to the case where $R$ places no mass at the origin, i.e.~$R>0$ a.s.

\begin{proposition} \label{prop:L1_properties}
Suppose $\boldsymbol X\eqdis R\boldsymbol S$ is $\ell_1$-norm symmetric with $R>0$ a.s. Then
\begin{itemize}
\item[(i)] $R\eqdis \|\boldsymbol X \|_1$ and $\boldsymbol S \eqdis \boldsymbol X / \|\boldsymbol X\|_1$
\item[(ii)] $\|\boldsymbol X \|_1$ and $\boldsymbol X / \|\boldsymbol X\|_1$ are independent.
\end{itemize}
\end{proposition}

\begin{proof}
The proof relies on the following fact: if $f:\R^n\to\R^m$ is a measurable map and $\boldsymbol Y \eqdis \boldsymbol Z$ are identically distributed $\R^n$-valued random variables, then $f(\boldsymbol Y) \eqdis f(\boldsymbol Z)$, see \cite{Fang/Kotz/Ng:1990}, p.~13. Letting $f:\boldsymbol y\mapsto (\|\boldsymbol y\|_1, \boldsymbol y / \|\boldsymbol y\|_1)$ for $\boldsymbol y\in \R^d\backslash\{\boldsymbol 0\}$ we get
$$
\left(\| \boldsymbol X\|_1, \frac{\boldsymbol X}{\|\boldsymbol X\|_1} \right) = f(\boldsymbol X) \eqdis f(R\boldsymbol S) = (R, \boldsymbol S),
$$
using that $\|\boldsymbol X\|_1>0$ a.s. Now $(i)$ clearly follows, and since we applied the result to the \emph{joint} distribution of $(R,\boldsymbol S)$, we also get $(ii)$. 
\proofend \end{proof}

The following characterization of $\ell_1$-norm symmetric distributions without point mass at the origin will be crucial in the sequel. Note that the proof of the direct half, $(i)$ implies $(ii)$, uses results on $d$-monotone functions. This fundamental connection between $d$-monotonicity and $\ell_1$-norm symmetric distributions is also evident when comparing the next result with Proposition \ref{prop:dmon_sf}.

\begin{proposition} \label{prop:L1_sf}
The following are equivalent:
\begin{itemize}
\item[(i)] $\bar H(\boldsymbol x)=f(x_1+\ldots+x_d)$, $\boldsymbol x\in \R^d_+$, is a survival function on the interior of $\R^d_+$.
\item[(ii)] $\bar H$ is the survival function of the $\ell_1$-norm symmetric random vector $\boldsymbol X \eqdis R\boldsymbol S$, $R>0$ a.s., and $f(x)=\P[RS_1>x]=\P[X_1>x]$.
\end{itemize}
\end{proposition}

\begin{proof}
$(ii) \Longrightarrow (i)$: Let $F$ be the distribution function of $R$. By Lemma \ref{lemma:condexpindep} we have
$$
\P[\boldsymbol X > \boldsymbol x] = \E\left[\P\left[\boldsymbol S > \frac{\boldsymbol x}{R} \Big| R\right]\right] = \int_{\R_+} \P\left[\boldsymbol S > \frac{\boldsymbol x}{t}\right] dF(t), \qquad \boldsymbol x\in\R^d_+,
$$
and similarly
$$
f(x) = \P[RS_1>x] =  \int_{\R^d_+} \P\left[S_1 > \frac{x}{t}\right] dF(t), \qquad x\geq0.
$$
It is well-known (see e.g.~Fang et al.~\cite{Fang/Kotz/Ng:1990}, Theorem 5.2 $(i)$) that the joint survival function of $\boldsymbol S\sim\mathrm{unif}(\mathcal S^d)$ is given by $\P[\boldsymbol S>\boldsymbol x] = (1-\|\boldsymbol x\|_1)^{d-1}_+$, $\boldsymbol x\in\R^d_+$. Hence $\P[S_1>x]=(1-x)^{d-1}_+$, $x\geq0$, and we get
$$
\P[\boldsymbol S>\boldsymbol x] = \P[S_1>\|\boldsymbol x\|_1], \qquad \boldsymbol x\in\R^d_+.
$$
This finally implies that
$$
\P[\boldsymbol X > \boldsymbol x] = \int_{\R_+} \P\left[S_1 > \frac{\|\boldsymbol x\|_1}{t}\right] dF(t) = f(\|\boldsymbol x\|_1), \qquad \boldsymbol x\in\R^d_+.
$$
$(i) \Longrightarrow (ii)$: Since $(i)$ holds by assumption, $f$ has to be $d$-monotone by Proposition \ref{prop:dmon_sf}. The representation theorem for $d$-monotone functions, Theorem \ref{thm:dmon_repr}, then gives a random variable $R$ with distribution function $F$ such that $f=\mathfrak W_d F$. Due to the invertibility of the Williamson $d$-transform (Proposition \ref{prop:Wd_prop}), $F$ is uniquely determined by $f$, and by the discussion after the proof of Theorem \ref{thm:dmon_repr}, $\mathfrak W_d F(x)=\P[RS_1>x]$ since we know that $S_1\sim\mathrm{Beta}(1,d-1)$. Finally, take $\boldsymbol S$ uniform on $\mathcal S^d$ and independent of $R$, and define $\boldsymbol X:=R\boldsymbol S$. As in the proof of the first part of the theorem it follows that $\bar{H}$ is the survival function of $\boldsymbol X$.
\proofend \end{proof}

\section{Archimedean copulas: characterization and stochastic representation} \label{section:Arch_char}

After the excursions in the two previous sections we now return to the Archimedean copulas and the problem of their characterization and representation. We will see how the theory of $d$-monotone functions and $\ell_1$-norm symmetric distributions enters the picture. The connection is actually little more than a simple observation; we use Sklar's theorem together with the fact that an Archimedean generator $\psi$ is a survival function on $\R_+$.

Suppose $C$ is an Archimedean copula generated by $\psi$. If $\psi$ is strict, the equality
$$
\bar H(\boldsymbol x) := C(\psi(x_1),\ldots,\psi(x_d)) = \psi(x_1+\ldots+x_d), \qquad \boldsymbol x\in\R^d_+
$$
certainly holds. If $\psi$ is non-strict, the equality holds if $x_i\leq\psi^{-1}(0)$ for all $i$, since in this case $\psi^{-1}\circ\psi(x_i)=x_i$. If on the other hand $x_i>\psi^{-1}(0)$ for some $i$, the right side is clearly zero, as is the left side since $\psi(x_i)=0$ and $C$ is a distribution function with uniform marginals (hence no mass on the boundary of $[0,1]^d$). So the equality holds also in this case. Sklar's theorem now implies that $\bar H$ is a joint survival function on $\R^d_+$, and even on the interior of $\R^d_+$ since $\psi(0)=1$ by the definition of Archimedean generators.

Conversely, if a $d$-dimensional survival function $\bar{H}$ on the interior of $\R^d_+$ has the form
$$
\bar H(\boldsymbol x) = \psi(x_1+\ldots+x_d), \qquad \boldsymbol x\in\R^d_+
$$
for some function $\psi$, the marginal survival functions are clearly $\psi(x)$, and Sklar's theorem gives the form of its survival copula $C$ as
$$
C(\boldsymbol u) = \psi(\psi^\leftarrow(u_1)+\ldots+\psi^\leftarrow(u_d)).
$$
If in addition $\psi$ is an Archimedean generator we have $\psi^\leftarrow=\psi^{-1}$, and $C$ is Archimedean by definition. We formulate these conclusions as

\begin{lemma} \label{lemma:char1}
The function $C:[0,1]^d\to [0,1]$ given by
\begin{equation} \label{eqn:cop1}
C(\boldsymbol u) = \psi(\psi^{-1}(u_1)+\ldots+\psi^{-1}(u_d)), \qquad \boldsymbol u\in[0,1]^d
\end{equation}
is a $d$-dimensional Archimedean copula if and only if
\begin{equation} \label{eqn:sf2}
\bar H(\boldsymbol x) = \psi(x_1+\ldots+x_d) = \psi(\|\boldsymbol x\|_1), \qquad \boldsymbol x \in \R^d_+
\end{equation}
is a joint survival function on the interior of $\R^d_+$. In this case, $C$ is the survival copula of $\bar H$.
\end{lemma}

\begin{proof}
That $C$ is the survival copula of $\bar H$ is immediate from the discussion above. The equivalence assertion also follows as soon as we prove that the assumption (\ref{eqn:sf2}) on $\bar H$ implies that $\psi$ is an Archimedean generator. Clearly $\psi(x)\to 0$ as $x\to\infty$ since $\bar H$ is a survival function. Moreover, $\psi(0)=1$ since $\bar H$ is a survival function on the interior of $\R^d_+$. To deduce continuity and strict monotonicity we must invoke Proposition \ref{prop:dmon_sf}, which says that $\psi$ is $d$-monotone. In particular, $\psi$ is convex and hence continuous on $(0,\infty)$, and since (\ref{eqn:sf2}) implies that $\psi$ is a marginal survival function, it must be right continuous everywhere which excludes a jump at $0$. Finally, since $\psi(x)\to0$ as $x\to\infty$, convexity implies that $\psi$ is strictly decreasing on the set where it is positive.
\proofend \end{proof}

It is clear from Lemma \ref{lemma:char1} that if $\psi$ generates a $d$-dimensional copula, it also generates $k$-dimensional copulas for all $k\leq d$. Indeed, $\psi(x_1+\ldots+x_k)$ is the $k$-dimensional marginal of $\bar H(\boldsymbol x)=\psi(x_1+\ldots+x_d)$, and we may apply Lemma \ref{lemma:char1} to this function.

We are now in a position to put the pieces together and give the main result that describes the structure of the Archimedean copulas. It is clear that the common denominator is the class of distributions whose survival functions can be written as a function of the $\ell_1$-norm of the argument. Due to Proposition \ref{prop:L1_sf}, these are precisely the $\ell_1$-norm symmetric distributions without a point mass at the origin.

\begin{theorem} \label{thm:main} (Main result)
The following are equivalent:
\begin{itemize}
\item[(i)] $C(\boldsymbol u)=\psi(\psi^{-1}(u_1)+\ldots+\psi^{-1}(u_d))$ is a $d$-dimensional Archimedean copula.
\item[(ii)] $\psi$ is a $d$-monotone Archimedean generator.
\item[(iii)] $\bar H(\boldsymbol x)=\psi(x_1+\ldots+x_d)$, $\boldsymbol x\in \R^d_+$, is a survival function on the interior of $\R^d_+$.
\item[(iv)] $\bar H$ is the survival function of the $\ell_1$-norm symmetric random vector $\boldsymbol X \eqdis R\boldsymbol S$, $R>0$ a.s., and $\psi(x)=\P[RS_1>x]=\P[X_1>x]$.
\end{itemize}
If $(i)$---$(iv)$ hold, $C$ is the survival copula corresponding to the survival function $\bar H$, and $\psi=\Wd F$, where $F$ is the distribution function of $R$.
\end{theorem}

\begin{proof}
The equivalence $(i)\Longleftrightarrow(iii)$ as well as the assertion that $C$ is the survival copula of $\bar H$ follow from Lemma \ref{lemma:char1}. Moreover, $(ii)\Longleftrightarrow(iii)$ and $(iii)\Longleftrightarrow(iv)$ follow from Proposition \ref{prop:dmon_sf} and Proposition \ref{prop:L1_sf}, respectively. Finally, since $S_1$ is $\mathrm{Beta}(1,d-1)$ distributed, $\psi=\P[RS_1>x]$ implies that
$\psi=\Wd F$ as in the discussion after the proof of Theorem \ref{thm:dmon_repr}.
\proofend \end{proof}

\begin{corollary} \label{cor:main}
Let $\boldsymbol U$ be distributed according to the $d$-dimensional Archimedean copula $C$ with generator $\psi$. Then $\boldsymbol X := (\psi^{-1}(U_1),\ldots,\psi^{-1}(U_d))$ is $\ell_1$-norm symmetric with radial part $R\sim F$ satisfying $R>0$ a.s.~and $\psi=\Wd F$.
\end{corollary}

\begin{proof}
The survival function of $\bar{H}$ of $\bX:=(\psi^{-1}(U_1),\ldots,\psi^{-1}(U_d))$ is
$$
\bar{H}(\bx) = \P[\psi^{-1}(U_1)>x_1,\ldots,\psi^{-1}(U_d)>x_d] = \P[U_1\leq \psi(x_1),\ldots,U_d\leq \psi(x_d)],
$$
since $\psi$ is decreasing. $\bU$ is distributed according to $C$, so
$$
\bar{H}(\bx) = \psi(\psi^{-1}(\psi(x_1))+\ldots+\psi^{-1}(\psi(x_1))) = \psi(\|\bx\|_1), \quad \bx\in\R^d_+,
$$
where the last equality follows as in the proof of Lemma \ref{lemma:char1}. The assertion now follows from Theorem \ref{thm:main}.
\proofend \end{proof}

Part $(ii)$ of Theorem \ref{thm:main} gives a necessary and sufficient condition that can be used to check whether a given Archimedean generator indeed generates a $d$-dimensional copula. One can for instance show that this is the case for the Gumbel generator $\psi^{\mathrm{Gu}}_\theta(x)$ in Example \ref{ex:arch_copulas} if $\theta\geq 1$. The Clayton generator  $\psi^{\mathrm{Cl}}_\theta(x)$ turns out to be $d$-monotone for $\theta \geq -1/(d-1)$, see McNeil and Ne\v{s}lehov\'a \cite{McNeil/Neslehova:2008}, Example 4. Part $(iv)$ gives a stochastic representation that greatly simplifies computations and can be used for simulation purposes. It is also the key to important interpretations of the Archimedean dependence structures, which is essential when justifying the use of these copulas in practice.

\begin{itemize}
\item[1.] We know that an $\ell_1$-norm symmetric random variable $\boldsymbol X\eqdis R\boldsymbol S$ has the properties that $\|\boldsymbol X\|_1\eqdis R$ and $\boldsymbol X/\|\boldsymbol X\|_1\eqdis \boldsymbol S$. Hence it can be viewed as an allocation of a random total amount $R$ of some resource among $d$ units, such that $S_i$ is the fraction allocated to unit $i$.
\item[2.] Moreover, $\|\boldsymbol X\|_1$ and $\boldsymbol X/\|\boldsymbol X\|_1$ are independent, meaning that the individual fractions do not depend on the total allocated amount.
\item[3.] Finally, requiring uniformity of $\boldsymbol S$ on the unit simplex $\mathcal S^d$ amounts to saying that all possible combinations of allocated fractions are equally likely to occur. This clearly imposes a strong symmetry on the distribution.
\end{itemize}

Archimedean copulas are the natural dependence structures of such ``fair-share'' distributions. One application from finance is in credit risk, where large portfolios of credit deals (e.g.~loans) are commonplace. In many cases such portfolios, or sub-portfolios, possess a large degree of homogeneity, which is consistent with point 3. Viewing the total loss as the ``resource'' (in credit risk it is common practice only to consider losses; $R$ will thus be positive) and assuming that the size of the loss does not influence its distribution between clients, we have 1.~and 2. The latter is reasonable only if the client exposures are of comparable orders of magnitude; knowing that the total loss was large, one could otherwise deduce that this large loss probably was caused by defaults among high-exposure clients rather than low-exposure clients. This would contradict independence between $\|\bX\|_1$ and $\bX/\|\bX\|_1$.

\section{A supplementary result} \label{section:Arch_compres}

Due to the invertibility of the Williamson $d$-transform, Theorem \ref{thm:main} implies that there is a one-to-one correspondence between a.s.~strictly positive random variables, or, equivalently,  $\ell_1$-norm symmetric distributions without a point mass at the origin, and $d$-monotone Archimedean generators. The correspondence between the generators and the Archimedean copulas is, however, \emph{not} one-to-one. Indeed, $\psi(x)$ and $\tilde \psi(x):=\psi(cx)$, $c>0$, give rise to the same copula. The corresponding radial distributions are $R=\mathfrak W_d^{-1}\psi$ and $\tilde R := R/c$. The next result shows that for a fixed copula $C$, the possible generators (or radial parts) cannot be more different than that.

\begin{proposition}
\begin{itemize}
\item[(i)] Two Archimedean generators $\psi$ and $\tilde\psi$ generate the same copula if and only if $\tilde\psi(cx)=\psi(x)$ for some $c>0$.
\item[(ii)] Let $\boldsymbol X$ and $\boldsymbol{\tilde X}$ be $\ell_1$-norm symmetrically distributed. They have identical survival copulas if and only if their radial parts $R$ and $\tilde R$ satisfy $\tilde R\eqdis R/c$ for some $c>0$.
\end{itemize}
\end{proposition}

\begin{proof}
Due to the one-to-one correspondence between the generators and radial parts as discussed above, it suffices to prove one of the two statements. We give a proof of $(i)$. Sufficiency follows immediately. For necessity, suppose $\psi$ and $\tilde \psi$ generate the same copula. Then
$$
\tilde \psi( \tilde\psi^{-1}(u_1) + \dots + \tilde\psi^{-1}(u_d)) = \psi( \psi^{-1}(u_1)+\dots+\psi^{-1}(u_d)), \qquad u_i\in[0,1].
$$
Now define $T:=\tilde \psi^{-1}\circ \psi$. Since the pseudo-inverse always is a right inverse, we have $T\circ\psi^{-1} = \tilde \psi^{-1}$ and $\tilde \psi \circ T = \psi$. The above equality is thus equivalent to
$$
\tilde \psi( T\circ\psi^{-1}(u_1) + \dots + T\circ\psi^{-1}(u_d)) = \tilde\psi\circ T( \psi^{-1}(u_1)+\dots+\psi^{-1}(u_d)),
$$
or, with $x_i:=\psi^{-1}(u_i)\in[0,x^*]$, where $x^*:=\psi^{-1}(0)$,
\begin{equation} \label{eqn:lspsi}
\tilde \psi( T(x_1) + \dots + T(x_d)) = \tilde\psi\circ T(x_1+\dots+x_d), \qquad x_i\in[0,x^*].
\end{equation}
We would now like to apply $\tilde\psi^{-1}$ on both sides of (\ref{eqn:lspsi}) to get a functional equation for $T$, but $\tilde\psi$ is only invertible on $\{\tilde\psi>0\}$. Now, if $x_1+\dots+x_d<x^*$ then $T(x_1+\dots+x_d)<\tilde\psi^{-1}(0)$, implying that the right side of (\ref{eqn:lspsi}) is strictly positive. The left side is then also strictly positive, and since $x_1+\dots+x_d<x^*$ implies that $x_i<x^*$ for each $i$ (all $x_i$ are positive) we conclude
$$
T(x_1) + \dots + T(x_d) = T(x_1+\dots+x_d), \qquad x_i>0, \ x_1+\dots+x_d<x^*.
$$
By continuity of $T$ this holds also for $x_i\geq 0$, $x_1+\dots+x_d< x^*$. With $x_1=\dots=x_d=0$ we get $dT(0)=T(0)$, so $T(0)=0$. Moreover, setting $x_3=\dots=x_d=0$ and using that $T(0)=0$ yields
$$
T(x_1) + T(x_2) = T(x_1+x_2), \qquad x_1,x_2\geq0, \ x_1+x_2< x^*.
$$
By Lemma \ref{lemma:fcneq} below, the only solutions to this functional equation are $T(x)=cx$, $x\in[0,x^*]$, where $c\in\R$. Hence $\psi(x)=\tilde\psi(cx)$ for $x\in[0,x^*]$, and since both $\psi$ and $\tilde \psi$ are Archimedean generators it is immediate that $c>0$. If $x^*=\infty$ we are done. Otherwise we note that $\psi(x^*)=0$, and so $\tilde\psi(x)=\psi(x/c)=0$ for $x\geq cx^*$, since $\tilde\psi$ is non-increasing.
\proofend \end{proof}

\begin{lemma} \label{lemma:fcneq}
Let $f$ be a real-valued, continuous function on $\R_+$ such that\footnote{Note the similarity to the well-known Cauchy functional equation: $f(x)+f(y)=f(x+y)$ for all $x\in \R$. See Bingham et al. \cite{Bingham/Goldie/Teugels:1987}, Section 1.1, for a treatment. However, although similar, our case requires a separate proof.}
\begin{equation} \label{eqn:fcneq2}
f(x)+f(y) = f(x+y), \qquad x,y\geq0, \ x+y< a,
\end{equation}
where $a\in(0,\infty]$. Then $f(x)=cx$ for all $x\in[0,a]$ and some $c\in\R$.
\end{lemma}

\begin{proof}
We prove the lemma in the case where $a<\infty$. If $a=\infty$ the proof is simpler in the sense that one does not have to be as careful to avoid going outside $[0,a]$, where the functional equation is defined. Otherwise it is very similar.

We may without loss of generality assume that $a=1$. Once this case is established, the general case follows by considering $\tilde f(x):=f(ax)$. Moreover, by continuity of $f$ we have (\ref{eqn:fcneq2}) for $x,y\geq0$, $x+y\leq a$. Now assume that
\begin{itemize}
\item[$(i)$] $f(nx) = nf(x)$ for all $x$ such that $0\leq nx\leq 1$
\end{itemize}
holds for some $n\in\mathbb N$, and take $x$ such that $0\leq (n+1)x\leq 1$. Then
$$
f(nx+x) = f(nx)+f(x) = nf(x)+f(x),
$$
using (\ref{eqn:fcneq2}) and $(i)$. The latter is trivially satisfied for $n=1$, so by induction it holds for all $n$. Take now $m\in\mathbb N$ and $x\in[0,1]$. Then $m\frac{x}{m}\in[0,1]$, so $(i)$ yields $f(x) = f(m\frac{x}{m}) = mf(\frac{x}{m})$, or, equivalently,
\begin{itemize}
\item[$(ii)$] $f(\frac{1}{m}x) = \frac{1}{m}f(x)$ for all $x\in[0,1]$.
\end{itemize}
Finally take $x\in[0,1]$ rational. Then $x=\frac{n}{m}$ for some $n\in\mathbb N \cup\{0\}$ and $m\in\mathbb N$ with $m\geq n$. Since $0\leq n\frac{1}{m}\leq 1$ and $\frac{1}{m}\leq 1$, we may apply $(i)$ and $(ii)$ to obtain
$$
f(x) = f\left(\frac{n}{m}\right) = nf\left(\frac{1}{m}\right) = \frac{n}{m}f(1)=xf(1).
$$
With $c:=f(1)$ we thus have $f(x)=cx$ for all \emph{rational} $x\in[0,1]$, and by density and continuity it also holds for all $x\in[0,1]$: take $x_n\to x\in[0,1]$ with $x_n\in[0,1]$ rational, and deduce that $f(x)=\lim_{n\to\infty} f(x_n)=cx$.
\proofend \end{proof}

\chapter{Regular variation and extreme value theory} \label{section:RV_EVT}

\section{Regular variation}

We summarize here some basic results from the theory of real-valued regularly varying functions. The main reference for this theory is Bingham et al.~\cite{Bingham/Goldie/Teugels:1987}, and the results are stated (and occasionally proven) in many other places as well. Two good references are Embrechts et al.~\cite{Embrechts/Kluppelberg/Mikosch:1997} and Resnick \cite{Resnick:2008}. The latter reference also discusses so-called multivariate regular variation in the context of multivariate extreme value theory.

We use the following notation: for two positive functions $f$ and $g$ we write $f(x)\sim g(x)$, $x\rightarrow\infty$, if $\lim_{x\rightarrow\infty}f(x)/g(x)=1$. Whenever $f(x)\sim 0$ appears, it is to be understood as $\lim_{x\to\infty}f(x)=0$. There will be no danger of confusing the two meanings of the symbol ``$\sim$''.

\begin{definition}
A positive, measurable function $f$ on $(0,\infty)$ is called regularly varying with index $\rho$, written $f\in \RV_\rho$, if for any $c>0$,
\[ \lim_{x\rightarrow\infty}\frac{f(cx)}{f(x)}=c^\rho. \]

It is called rapidly varying (with index $-\infty$), $f\in \RV_{-\infty}$, if
\[ \lim_{x\rightarrow\infty}\frac{f(cx)}{f(x)} = \left\{ \begin{array}{ccc}
    0 & if \ \ c>1   \\
    \infty & if \ \ 0<c<1.  
\end{array} \right. \]

Functions in $\RV_0$ are called slowly varying.
\end{definition}

Note that a regularly varying function $f\in\RV_\rho$ always can be represented as $f(x)=x^\rho L(x)$, with $L\in \RV_0$; just set $L(x):=f(x)x^{-\rho}$. In this thesis, the function $f$ will usually be the tail of some distribution function $F$, i.e.~$f=1-F=\bar F$. We will then sometimes write $X\in \RV_\rho$ for a random variable $X\sim F$ to denote that $\bar F\in \RV_\rho$.

\begin{theorem} \label{thm:Karamata}
(Karamata's theorem) Let $L\in RV_0$ be locally bounded in $[0,\infty]$. Then
\begin{itemize}
\item[(i)] If $\sigma<0$, then $\int_x^\infty t^{\sigma-1} L(t)dt \sim -x^{\sigma}L(x)/\sigma$, $x\rightarrow\infty$.
\item[(ii)] If $\sigma>0$, then $\int_0^x t^{\sigma-1} L(t)dt \sim x^{\sigma}f(x)/\sigma$, $x\rightarrow\infty$.
\item[(iii)] The function $x\mapsto \int_0^x t^{-1} L(t)dt$ is slowly varying, if finite.
\end{itemize}
\end{theorem}

\begin{proof}
These are Theorem 1.5.10, 1.5.8 and 1.5.9a in Bingham et al.~\cite{Bingham/Goldie/Teugels:1987}. 
\proofend \end{proof}

\begin{theorem} (Monotone Density Theorem) \label{thm:mon_density}
Suppose $U(x)\sim cx^\rho L(x)$, $x\rightarrow\infty$, where $c\in\R$, $\rho\in\R$, $L\in RV_0$ and suppose $u$ is some eventually monotone function.
\begin{itemize}
\item[(i)] If $U(x)=\int_x^\infty u(y)dy$, then $u(x)\sim -c\rho x^{\rho-1}L(x)$ as $x\longrightarrow\infty$.
\item[(ii)] If $U(x)=\int_0^x u(y)dy$, then $u(x)\sim c\rho x^{\rho-1}L(x)$ as $x\longrightarrow\infty$.
\end{itemize}
\end{theorem}

\begin{proof}
\emph{(ii)} is Theorem 1.7.2 in Bingham et al.~\cite{Bingham/Goldie/Teugels:1987}, and the proof of \emph{(i)} is virtually identical. 
\proofend \end{proof}

\begin{theorem} \label{thm:rapvar}
Suppose $f\in RV_{-\infty}$ is non-increasing and take $\sigma\in\R$. Then $\int_{x_0}^\infty t^\sigma f(t)dt<\infty$ for some $x_0\geq0$, and
\[ \lim_{x\rightarrow\infty} \frac{\int_x^\infty t^\sigma f(t)dt}{x^{\alpha+1}f(x)} = 0. \]
\end{theorem}

\begin{proof}
This is (part of) Theorem A3.12 in Embrechts et al.~\cite{Embrechts/Kluppelberg/Mikosch:1997}.
\proofend \end{proof}

\section{Extreme value theory}

Univariate extreme value theory in its most basic form is concerned with the study of limit relations of the form
\begin{equation} \label{eqn:evt1}
F^n(c_n x + d_n) \to G(x) \ \  \mathrm{weakly}, \qquad n\to\infty,
\end{equation}
where $F$ and $G$ are distribution functions with $G$ non-degenerate, $c_n>0$, and $d_n\in \R$. The probabilistic interpretation is the following: Suppose $X_1,X_2,\ldots$ are i.i.d.~random variables distributed according to $F$, and let $M_n:=\max(X_1,\ldots,X_n)$ denote the largest of the first $n$ variables in the sequence. Then, by independence, $\P[M_n\leq x] = \P^n[X_1\leq x]=F^n(x)$, and after translation and scaling by $d_n$ and $c_n$, $\P[c_n^{-1}(M_n-d_n)\leq x] = F^n(c_nx+d_n)$. The relation (\ref{eqn:evt1}) thus expresses \emph{weak convergence of the normalized maxima} $c_n^{-1}(M_n-d_n)$.

If (\ref{eqn:evt1}) holds for some specific $F$ and non-degenerate $G$, we say that $F$ lies in the \emph{maximum domain of attraction} (or just \emph{domain of attraction}) of $G$, denoted $F\in\MDA(G)$. The first issue to resolve is whether the limit $G$ is unique. More precisely, is it possible to find normalizing sequences in (\ref{eqn:evt1}) so that for the same $F$ we get convergence to two different limits $G_1$ and $G_2$? It turns out that this can happen, but in that case $G_1$ and $G_2$ are of the same \emph{type}, meaning that $G_1(x)=G_2(Ax+B)$ for some constants $A>0$ and $B\in\R$. This is part of the Convergence to Types Theorem (see Proposition 0.2 in Resnick \cite{Resnick:2008}).

The next task is now to determine the possible limits $G$. This was done by Fisher and Tippett \cite{Fisher/Tippett:1928}, and it turns out that the limit (up to type) belongs to one of three families:
\begin{itemize}
\item[$(i)$] The Fr\'echet family: $\Phi_\alpha(x)=\exp\{-x^{-\alpha}\}$, $x\geq0$, for some $\alpha\in(0,\infty)$.
\item[$(ii)$] The Gumbel distribution: $\Lambda(x) = \exp\{-e^{-x}\}$, $x\in\R$.
\item[$(iii)$] The Weibull family: $\Psi_\alpha(x)=\exp\{-(-x)^\alpha\}$, $x<0$, for some $\alpha\in (0,\infty)$.
\end{itemize}

For each of these three families it now remains to characterize the respective domains of attraction and the corresponding normalizing sequences $(c_n)$ and $(d_n)$.

\begin{theorem} \label{thm:MDA_char}
Let $F$ be a distribution function.
\begin{itemize}
\item[(i)] (Fr\'echet case) $F\in\MDA(\Phi_\alpha)$ if and only if $\bar{F}\in RV_{-\alpha}$ for $\alpha\in(0,\infty)$.
\item[(ii)] (Gumbel case) $F\in\MDA(\Lambda)$, if and only if there exists a positive function $a$, called the \emph{auxiliary function} of $F$, with $a(x)/x\to 0$ as $x\to x^*$, where $x^*$ is the right endpoint of $F$, such that
\begin{equation} \label{eqn:gamma}
\lim_{x\uparrow x^*} \frac{\bar F(x+a(x)t)}{\bar F(x)} = e^{-t}, \qquad \forall t\in\R.
\end{equation}
In this case one may choose $a(x) := \int_x^{x^*} \bar F(t) dt / \bar F(x)$. If $F$ has a density, we may alternatively choose $a(x):= -\bar F(x)/\bar F'(x)$. It is always possible to find an auxiliary function that is absolutely continuous with $\lim_{x\to x^*} a'(x) = 0$.
\item[(iii)] (Weibull case) $F\in\MDA(\Psi_\alpha)$ if and only if $F$ has finite right endpoint $x^*$ and $\bar{F}(x^*-x^{-1})\in RV_{-\alpha}$ for $\alpha\in(0,\infty)$.
\end{itemize}
\end{theorem}

\begin{proof}
For $(i)$ and $(iii)$, see Gnedenko \cite{Gnedenko:1943}, Theorem 4 and Theorem 5, respectively. For the first part of $(ii)$, see Gnedenko \cite{Gnedenko:1943}, Theorem 7. For the two choices of auxiliary function, see de Haan \cite{deHaan:1971}, Lemma 5 and Resnick \cite{Resnick:2008}, Proposition 1.1 (b). That the auxiliary function can be taken absolutely continuous with the given property can be deduced from Resnick \cite{Resnick:2008}, Proposition 1.4 and the proof of Proposition 1.1(a) in the same reference.
\proofend \end{proof}

We note that in the Gumbel case, all auxiliary functions of a given distribution $F$ are tail equivalent. This convenient fact is proved in \cite{Resnick:2008}, p.~26\footnote{Alternatively one may invert both sides of (\ref{eqn:gamma}) to obtain $\{\bar F^{\leftarrow}(\bar F(x)t)-x\}/a(x)\to -\log t$ as $x\uparrow x^*$. If this also holds with $a(x)$ replaced by some $\hat a(x)$, then clearly $a(x)/\hat a(x)\to 1$. That inverting (\ref{eqn:gamma}) is allowed can be proven in the same way as  Proposition 0.1 in Resnick \cite{Resnick:2008}.}. Moreover, if $F\in\MDA(\Lambda)$ has unbounded support, then $\bar F$ is rapidly varying, $\bar F\in\RV_{-\infty}$.

In this thesis we will not need the full answer to the question of how to find appropriate normalizing sequences. For this we refer to Resnick \cite{Resnick:2008}, Proposition 1.11 and 1.13 for the Fr\'echet and Weibull domains, respectively, and Proposition 1.1, 1.4 and 1.19 for the Gumbel domain. For our purposes the following is sufficient:

\begin{proposition} \label{prop:normseq}
If $F\in\MDA(\Psi_\alpha)$,  the normalizations $c_n=x^*-(1/\bar F)^{\leftarrow}(n)$ and $d_n=x^*$ yield $G=\Psi_\alpha$ in (\ref{eqn:evt1}). Here $x^*$ is the right endpoint of $F$.
\end{proposition}

\begin{proof}
This follows from Proposition 1.13 in Resnick \cite{Resnick:2008}, or can be deduced from the proof of Theorem 5 in Gnedenko's paper \cite{Gnedenko:1943}.
\proofend \end{proof}

\chapter{Tail properties of the Williamson $d$-transform} \label{section:Wd_tail}

We are interested in the extremal behavior of Archimedean copulas, which depends on the generator near its right endpoint and near the origin. We prefer, however, to work with $\ell_1$-norm symmetric distributions instead of Archimedean generators, and it is thus of interest to connect the extremal behavior of the radial distributions to that of the corresponding generators. By Theorem \ref{thm:main} this amounts to investigating the tail properties of the Williamson $d$-transform. This program is carried out in the present chapter. We separate the analysis into four different cases: regular variation at infinity, the Gumbel domain of attraction, the behavior near the right endpoint in the case of bounded support, and the behavior at the origin.

Due to the interpretation of the Williamson $d$-transform as the survival function of certain Beta mixtures, the proofs of some of the results below can be significantly simplified using the following result by Breiman, slightly extended to fit our needs:

\begin{lemma} \label{lemma:Breiman}
Let $X$ and $Y$ be two non-negative random variables such that $\bar F \in \RV_{-\alpha}$, $\alpha\in(0,\infty)$, where $F$ is the distribution function of $X$, and $\E[Y^{\alpha+\varepsilon}]<\infty$ for some $\varepsilon>0$. Then
$$
\frac{\P[XY>x]}{\P[X>x]} \to \E[Y^\alpha], \qquad x\to \infty.
$$
If $Y\in(0,1)$ a.s., the conclusion holds also with $\alpha=0$ and $\alpha=\infty$, where in the latter case $\E[Y^\alpha]=0$ by convention.
\end{lemma}

\begin{proof}
The case when $\alpha\in(0,\infty)$ is covered by Breiman's original result \cite{Breiman:1965}, so we suppose that $Y\in (0,1)$ a.s. Denoting by $G$ the distribution function of $Y$ we get, using Lemma \ref{lemma:condexpindep},
$$
\lim_{x\to\infty} \frac{\P[XY>x]}{\P[X>x]} = \lim_{x\to\infty} \int_0^1 \frac{\bar F(x/t)}{\bar F(x)} dG(t).
$$
Since $\bar F$ is non-increasing and $t\in(0,1)$, the integrand is bounded by one, so we may interchange limit and integration by the Dominated Convergence Theorem. If $\alpha=0$, $\bar F(x/t)/\bar F(x) \to 1$ as $x\to\infty$. If $\alpha=\infty$, the limit is zero. In either case we clearly reach the desired conclusion.
\proofend \end{proof}

As an aid for the intuition, the reader is encouraged to think of the results in this chapter as probabilistic statements about mixture variables, and not only as analytical properties of integral transforms.

\section{Regular variation to the right, strict case} \label{section:regvar_strict}

A strict Archimedean generator $\psi$, viewed as a survival function, can be expected to possess interesting extremal behavior only if $\psi\in \RV_{-\alpha}$ for some $\alpha\in(0,\infty]$; this follows from Theorem \ref{thm:MDA_char}, since any other behavior of $\psi$ would exclude it from being in some domain of attraction. By definition, $\psi\in \RV_{-\alpha}$ means that

\begin{equation} \label{eqn:lim_strict}
\lim_{x\to\infty}\frac{\psi(cx)}{\psi(x)} = c^{-\alpha}, \qquad \forall c>0.
\end{equation}

Now, by Theorem \ref{thm:main}, $\psi$ has to be $d$-monotone in order to generate a copula in dimension $d$. At first sight, this seems to be a fairly strong regularity condition on $\psi$, especially for large $d$, and one wonders if perhaps (\ref{eqn:lim_strict}) always holds, at least if we also allow $\alpha=0$. This question is particularly relevant in view of the fact that all commonly encountered families of Archimedean copulas have generators that do satisfy (\ref{eqn:lim_strict}). We now give a counterexample to show that this is not always the case.

\begin{example} \label{ex:no_MDA}
Define
$$
\psi(x) := \frac{1 + a\sin(\log(1+x))}{1+x}, \qquad x\geq 0,
$$
for some $a>0$ to be determined. Then $\psi(0)=1$ and $\lim_{x\to\infty}\psi(x)=0$, and if $a<1$, then $\psi(x)>0$ for all $x\geq0$. We show that for any $d$, one can choose $a>0$ such that $\psi$ is $d$-monotone. Suppose that for some $k\in\{0,1,\dots\}$,
\begin{equation} \label{eqn:ex2}
\psi^{(k)}(x) = (-1)^k \frac{k! + a c_k \sin(\varphi_k + \log(1+x))}{(1+x)^{k+1}}, \qquad x\geq 0,
\end{equation}
where $c_k$ and $\varphi_k$ are constants. Differentiating and using standard trigonometric identities one obtains
\begin{eqnarray*}
\psi^{(k+1)}(x) &=& (-1)^{k+1} \frac{(k+1)! + a(k+1)c_k\sin(\varphi_k + \log(1+x)) - a c_k \cos(\varphi_k + \log(1+x))}{(1+x)^{k+2}} \\
&=& (-1)^{k+1} \frac{(k+1)! + a c_{k+1} \sin(\varphi_{k+1} + \log(1+x))}{(1+x)^{k+2}},
\end{eqnarray*}
where $c_{k+1} = c_k \sqrt{k^2+2k+2}$ and $\varphi_{k+1}=\varphi_k + \arctan(1/(k+1))$. Since (\ref{eqn:ex2}) holds for $k=0$ with $c_0=1$ and $\varphi_0=0$, it follows by induction that it holds for all $k\in\{0,1,\dots\}$. Since $\sin(y)\geq -1$ for all $y\in\R$, it is clear that given $d\geq 2$, $(-1)^k\psi^{(k)}(x)\geq 0$ for all $k=1,\dots,d-1$, provided $0<a \leq k!/c_{k}$ for each $k$.  Since for $k\geq 2$ it holds that $c_k>k!$ and $\max(c_0,\dots,c_k)=c_k$, such an $a$ can indeed be found and satisfies $a<1$ (take $a=(d-1)!/c_{d-1}$), implying that $\psi$ is indeed $d$-monotone.

We now consider the behavior of $\psi(cx)/\psi(x)$ as $x\to\infty$ along the sequence $x_n=e^{n\pi}$. We write
$$
\frac{\psi(cx_n)}{\psi(x_n)} = \frac{1 + a\sin(\log(1+cx_n))}{1 + a\sin(\log(1+x_n))} \frac{1+x_n}{1+cx_n} =: A(x_n) \frac{1+x_n}{1+cx_n}.
$$
Since $(1+x_n)/(1+cx_n)\to 1/c$, we focus attention on $A(x_n)$. With $x_n$ as above we get
$$
A(x_n) =  \frac{1 + a\sin(n\pi + \log(e^{-n\pi}+c))}{1 + a\sin(n\pi + \log(e^{-n\pi}+1))} =  \frac{1 + a(-1)^n\sin(\log(e^{-n\pi}+c))}{1 + a(-1)^n\sin(\log(e^{-n\pi}+1))}.
$$
The denominator converges to one as $n\to\infty$, but the numerator does not converge in general: along the subsequence where all $n$ are even, it converges to $1+a\sin(\log c)$, while if all $n$ are odd, the limit is $1-a\sin(\log c)$. So unless $\log c = m\pi$ for some integer $m$, $\psi(cx_n)/\psi(x_n)$ does not converge.
\end{example}

We remark that there are known families of Archimedean copulas where the generators satisfy (\ref{eqn:lim_strict}) with $\alpha=0$. One such example is $\psi(x)=\theta/\log(e^\theta + x)$ with $\theta\in (0,\infty)$, which generates bivariate Archimedean copulas (see Nelsen \cite{Nelsen:2006}, item 19 of Table 4.1, where however the inverse generator $\phi(t)=e^{\theta/t}-e^\theta$ is given instead of $\psi$). In such cases it is still possible to draw conclusions about the tail behavior of the copula, see Charpentiers and Segers \cite{Charpentier/Segers:2008}, although the corresponding $\ell_1$-norm symmetric distributions do not possess non-degenerate extremal behavior (under affine normalizations).

After this introductory discussion, we turn to sufficient and in most cases necessary conditions on the radial part of the corresponding $\ell_1$-norm symmetric distribution for (\ref{eqn:lim_strict}) to hold.

\begin{proposition} \label{prop:Wd_tail_s1}
Let $F$ be a distribution function with $F(0)=0$ and let $\psi=\Wd F$, $d\geq 2$, be its Williamson $d$-transform. Then
\begin{itemize}
\item[(i)] If $\bar F\in RV_{-\alpha}$, $\alpha\in[0,\infty)$, then $\psi\in\RV_{-\alpha}$.
\item[(ii)] If $\psi\in RV_{-\alpha}$, $\alpha\in(0,\infty)$, then $\bar F\in\RV_{-\alpha}$.
\end{itemize}
\end{proposition}

\begin{proof}
Suppose $R\sim F$ and let $S_1$ be $\mathrm{Beta}(1,d-1)$ distributed and independent of $R$. Then, as we have seen (e.g.~in Theorem \ref{thm:main}), $\psi(x)=\P[RS_1>x]$. Observing that  $S_1\in(0,1)$ a.s., we obtain from Lemma \ref{lemma:Breiman} that $\bar F(x)\sim B_1 \psi(x)$ as $x\to\infty$ for some $B_1>0$, whence the statement $(i)$.

For part $(ii)$, we have by assumption that $\psi(x)=x^{-\alpha}L(x)$ for some $\alpha>0$ and slowly varying $L$. We first show by induction that, for $k=0,1,\ldots,d-1$,
$$
(-1)^kx^k\psi^{(k)}(x) \sim c_k x^{-\alpha}L(x) = c_k\psi(x), \qquad x\to\infty,
$$
where $c_0=1$ and $c_k=\prod_{j=0}^{k-1}(j+\alpha)$ for $k=1,\dots,d-1$. For $k=0$ the statement is trivially true. For $k=0,\ldots,d-2$,
$$
\psi^{(k)}(x) = -\int_x^\infty \psi^{(k+1)}(y)dy,
$$
since $\psi^{(k)}(x)\rightarrow0$ as $x\rightarrow \infty$. Moreover,  $\psi^{(k+1)}$ is monotone since $\psi$ is $d$-monotone, and $\psi^{(k)}(x)\sim (-1)^{k}c_{k}x^{-k-\alpha}L(x)$ by the induction assumption. The Monotone Density Theorem now implies that $\psi^{(k+1)}(x)\sim (-1)^kc_k(-k-\alpha)x^{-k+1-\alpha}L(x)$, which is precisely the statement we wanted to prove.

Inserting this into the inverse Williamson $d$-transform yields
$$
\frac{\bar F (x)}{\psi(x)} = \frac{1-\Wd^{-1}\psi(x)}{\psi(x)} = \sum_{k=0}^{d-1}\frac{1}{k!}\frac{(-1)^kx^k\psi^{(k)}(x)}{\psi(x)} \to \sum_{k=0}^{d-1} \frac{1}{k!} c_k, \quad x\to\infty.
$$
The limit is clearly positive since $\alpha>0$.
\proofend \end{proof}

\textbf{Remark.} One can show that  
$$
\sum_{k=0}^{d-1} \frac{1}{k!} c_k = 1 +  \sum_{k=1}^{d-1} \frac{1}{k!} \prod_{j=0}^{k-1}(j+\alpha) = \prod_{k=1}^{d-1}\left(1+\frac{\alpha}{k}\right),
$$
which by Breiman's result is equal to $1/\E[S_1^\alpha]$. As a by-product we thus get an expression for the raw moments of $\mathrm{Beta}(1,d-1)$ distributed random variables.

\section{The Gumbel domain of attraction} \label{section:rapvar}

If $\bar F\in \RV_{-\infty}$, which for instance is the case if $F\in \MDA(\Lambda)$ and has unbounded support, we may use Lemma \ref{lemma:Breiman} similarly as in the proof of Proposition \ref{prop:Wd_tail_s1} $(i)$ to conclude that $\psi(x)/\bar F(x)\to 0$ as $x\to \infty$. This gives no information about the precise behavior of the tail of $\psi$. In what follows we derive more precise results in the case where $F$ and/or $\psi$ is in the Gumbel domain of attraction. This is the relevant case in the context of extreme value theory. It turns out to be unnecessary to distinguish between the strict and the non-strict cases.

We proceed by proving a number of intermediate results. The condition (\ref{eqn:gamma}) from Theorem \ref{thm:MDA_char} will appear several times, and we emphasize that the requirements on the auxiliary function $a$ are part of this condition. Namely, $a(x)>0$ and $a(x)/x\to 0$ as $x$ tends to the right endpoint $x^*$ (which may be infinity).

\begin{lemma} \label{lemma:gamma_diff}
Suppose the non-negative, non-increasing function $f$ satisfies (\ref{eqn:gamma}), i.e.
$$
\lim_{x\uparrow x^*} \frac{f(x+a(x)t)}{f(x)} = e^{-t}, \qquad \forall t\in\R,
$$
where $a(x)>0$ and $a(x)/x\to 0$ as $x\uparrow x^*$. If $f$ has a density $f'$ that is monotone in a neighborhood of $x^*$, then
$$
\lim_{x\uparrow x^*} \frac{f'(x+a(x)t)}{f'(x)} = e^{-t}, \qquad \forall t\in\R,
$$
\end{lemma}

\begin{proof}
The core of the argument is similar to the proof of Theorem 3.10.11 in Bingham et al.~\cite{Bingham/Goldie/Teugels:1987}. For all $x$ large enough, $x+a(x)t = x(1+ta(x)/x)$ is close to $x^*$ for any $t\in\R$, since $a(x)/x\to 0$ as $x\uparrow x^*$. Now choose any $u,v\in\R$, $u>v$. For all $x$ sufficiently close to $x^*$, say $x\in(x_0,x^*)$, we thus have that $x+a(x)v$ and $x+a(x)u$ are both in the region where $f'$ is monotone. Since $f$ is non-increasing, $-f'$ is non-negative. We look at the case where $-f'$ is non-increasing near $x^*$, the non-decreasing case being completely analogous. We obtain
\begin{eqnarray*}
-\frac{f(x+a(x)u)-f(x+a(x)v)}{f(x)} &=& \frac{a(x)}{f(x)}\int_v^u -f'(x+a(x)t)dt \\
&\leq& -\frac{a(x)}{f(x)}f'(x+a(x)v)(u-v),
\end{eqnarray*}
where each expression is non-negative. The left side tends to $e^{-v}-e^{-u}$ as $x\to x^*$, and hence
$$
\liminf_{x\to x^*}\ -\frac{a(x)}{f(x)}f'(x+a(x)v) \geq \lim_{u\downarrow v} \frac{e^{-v}-e^{-u}}{u-v} = e^{-v},
$$
where the last equality follows by e.g.~l'H\^opital's rule. Similarly one shows
$$
\limsup_{x\to x^*} \ -\frac{a(x)}{f(x)}f'(x+a(x)u) \leq e^{-u},
$$
and so
\begin{equation} \label{eqn:difflema1}
\lim_{x\to x^*} \ -\frac{a(x)}{f(x)}f'(x+a(x)t) = e^{-t}, \qquad \forall t\in\R.
\end{equation}
In particular, we may take $t=0$ to obtain $-f(x)/a(x)\sim f'(x)$, $x\to x^*$, which implies that $f'$ is non-zero near $x^*$. Substituting back into (\ref{eqn:difflema1}) gives the conclusion.
\proofend \end{proof}

\begin{lemma} \label{lemma:gamma_pow}
Suppose the function $f$ satisfies (\ref{eqn:gamma}). Then, for any $\rho\in\R$, $g(x)=x^\rho f(x)$ also satisfies (\ref{eqn:gamma}) with the same auxiliary function.
\end{lemma}

\begin{proof}
A simple computation yields
$$
\frac{g(x+a(x)t)}{g(x)} = \left(1+\frac{a(x)}{x}t\right)^\rho \frac{f(x+a(x)t)}{f(x)} \to e^{-t},
$$
since $f$ satisfies (\ref{eqn:gamma}), which in particular means that $a(x)/x\to 0$ as $x\to x^*$.
\proofend \end{proof}

\begin{lemma} \label{lemma:gamma_int}
Suppose the non-negative, non-increasing function $f$ satisfies (\ref{eqn:gamma}). Then $g(x)=\int_x^\infty f(t)dt$ also satisfies (\ref{eqn:gamma}) with the same auxiliary function.
\end{lemma}

\begin{proof}
Suppose $\hat a$ is the auxiliary function of $f$. By Theorem \ref{thm:MDA_char} we can find another auxiliary function $a$ that is absolutely continuous with $a'(x)\to 0$, $x\to x^*$. L'H\^opital's rule then gives
$$
\lim_{x\uparrow x^*} \frac{g(x+a(x)t)}{g(x)} = \lim_{x\uparrow x^*} \frac{-f(x+a(x)t)}{-f(x)}(1+a'(x)t).
$$
The limit on the right is clearly $e^{-t}$ by the properties of $f$ and $a$. So $g$ satisfies (\ref{eqn:gamma}) with auxiliary function $a$, hence also with $\hat a$, since the two are tail equivalent, i.e.~$\hat a(x)/a(x)\to 1$, $x\to x^*$.
\proofend \end{proof}

\begin{lemma} \label{lemma:psiFprim}
Let $F$ be a distribution function with $F(0)=0$ and let $\psi = \Wd F$. Then $\psi^{(d-1)}$ is differentiable if and only if $F$ is differentiable, and in this case
$$
\bar F'(x) = \frac{(-1)^{d-1}}{(d-1)!}x^{d-1}\psi^{(d)}(x).
$$
\end{lemma}

\begin{proof}
From the inversion formula of the Williamson $d$-transform in Proposition \ref{prop:Wd_prop}, it is clear that  $\psi^{(d-1)}$ is differentiable if and only if $F$ is, and that the derivative $\bar F'$ in this case is given by $\bar F'(x)=f_d(x)$, where we define
$$
f_n(x) := D \left[ \sum_{k=0}^{n-1} \frac{(-1)^k}{k!}x^k \psi^{(k)}(x) \right]
$$
for all $n$ such that $\psi^{(n)}$ exists. Here $D[\cdot]$ denotes differentiation. We are done if we can show that for any such $n$,
\begin{equation} \label{eqn:psiFprim1}
f_n(x) = \frac{(-1)^{n-1}}{(n-1)!}x^{n-1}\psi^{(n)}(x).
\end{equation}
We argue by induction over $n$. For $n=2$, $f_n(x)=D[\psi(x)-x\psi'(x)] = -x\psi''(x)$, so (\ref{eqn:psiFprim1}) holds. Suppose (\ref{eqn:psiFprim1}) holds for $n-1$ and $\psi^{(n)}$ exists. Then
\begin{eqnarray*}
f_n(x) &=& f_{n-1}(x) + D\left[ \frac{(-1)^{n-1}}{(n-1)!} x^{n-1}\psi^{(n-1)}(x) \right] \\
&=&  \frac{(-1)^{n-1}}{(n-1)!} x^{n-1}\psi^{(n)}(x).
\end{eqnarray*}
\proofend \end{proof}

We are now ready to prove the following proposition concerning the connection between $\bar F$ and $\psi$ when either is in the Gumbel domain of attraction. Note the assumption that $\bar F$ or $\psi$ be differentiable---we are not treating the most general case here.

\begin{proposition} \label{prop:Wd_tail_Gumbel}
Let $F$ be a distribution function with $F(0)=0$, let $\psi=\Wd F$ for $d\geq 2$ be its Williamson $d$-transform, and suppose that either $F$ or $\psi^{(d-1)}$ is differentiable. Both are then differentiable, and if we in addition assume that $F'$ is monotone in a neighborhood of its right endpoint $x^*$, then $F\in\MDA(\Lambda)$ if and only if $\psi\in\MDA(\Lambda)$.
In this case,
$$
a(x) = \frac{\int_x^{x^*}\bar F(t)dt}{\bar F(x)} \qquad \mathrm{and} \qquad a(x) = \frac{\int_x^{x^*} \psi(t)dt}{\psi(x)}
$$
are possible choices of auxiliary function for both $\bar F$ and $\psi$.
\end{proposition}

\begin{proof}
Lemma \ref{lemma:psiFprim} implies that both $\bar F$ and $\psi$ are differentiable and satisfy
\begin{equation} \label{eqn:WdtailG}
\bar F'(x) = \frac{(-1)^{d-1}}{(d-1)!}x^{d-1}\psi^{(d)}(x).
\end{equation}
We start with necessity, so suppose $\bar F$ is in the Gumbel domain. Theorem \ref{thm:MDA_char} and the monotonicity assumption on $F'=-\bar F'$ then enable us to use Lemma \ref{lemma:gamma_diff} to conclude that $-\bar F'$, and hence $(-1)^d x^{d-1}\psi^{(d)}$, satisfies (\ref{eqn:gamma}). The auxiliary function carries over, and Lemma \ref{lemma:gamma_pow} implies that $(-1)^d\psi^{(d)}$ also satisfies (\ref{eqn:gamma}) with the same auxiliary function. Since
$$
(-1)^k \psi^{(k)}(x) = \int_x^\infty (-1)^{k+1} \psi^{(k+1)}(t)dt
$$
for all $k=0,1,\dots,d-1$, repeated use of Lemma \ref{lemma:gamma_int} yields that $\psi$ satisfies (\ref{eqn:gamma}), again with the same auxiliary function. Hence $\psi\in\MDA(\Lambda)$.

To prove sufficiency the procedure is simply reversed. In this direction the use of Lemma \ref{lemma:gamma_diff} is justified since $\psi^{(k)}$ are monotone by $d$-monotonicity for $k=1,\ldots,d-2$, and also for $k=d-1$ since $\psi^{(d-2)}$ is convex. For $k=d$, we use (\ref{eqn:WdtailG}) together with the eventual monotonicity of $\bar F'$, which implies that $\psi^{(d)}$ also is eventually monotone---indeed, suppose $f(x)=x^\rho g(x)$ for some functions $f$ and $g$ on $(0,x^*)$, where $\rho\geq 0$. If for any $x_0<x^*$ one can find $x,y\in(x_0,x^*)$, $x<y$, such that $g(x)<g(y)$, then $f$ cannot be non-increasing on any left neighborhood of $x^*$, since clearly $f(x)<f(y)$. With $f= -\bar F'$ and $g=-\psi^{(d)}$, this gives a contradiction if $\psi^{(d)}$ is not eventually monotone.
\proofend \end{proof}

\section{Regular variation to the right, non-strict case}

We start with a result that describes how non-strict $d$-monotone Archimedean generators and their derivatives behave in case they are regularly varying near the right endpoint.

\begin{proposition} \label{prop:nonstrict_regvar}
Let $\psi=\Wd F$ for some distribution function $F$ with $F(0)=0$, and suppose the right endpoint $x^*=\psi^{-1}(0)$ is finite. 
\begin{itemize}
\item[(i)] If $\psi(x^*-x^{-1})\in \RV_{-\alpha}$, $\alpha>0$, then
$$
(-1)^k \psi^{(k)}(x^*-x^{-1}) \in RV_{-\alpha+k}
$$
for $k=1,\ldots,d-1$. Moreover, in this case $\alpha\geq d-1$.
\item[(ii)] If $(-1)^{d-1}\psi^{(d-1)}(x^*-x^{-1}) \in \RV_{-\alpha}$, $\alpha \geq 0$, then
$$
(-1)^k \psi^{(k)}(x^*-x^{-1}) \in \RV_{-\alpha-(d-1)+k}
$$
for $k=0,\dots,d-2$.
\end{itemize}
\end{proposition}

\begin{proof}
We first prove $(i)$. Since $\psi^{(k)}$ is continuous for $k=0,\ldots,d-2$, $\psi^{(k)}(x^*)=0$ for these $k$. Hence for $k=1,\ldots,d-1$,
\begin{equation} \label{eqn:nonstr_regv}
\psi^{(k-1)}(x^*-x^{-1}) = -\int_x^\infty t^{-2}\psi^{(k)}(x^*-t^{-1})dt.
\end{equation}
Note that $\psi^{(d-1)}$ is in general only a density, not a derivative, of $\psi^{(d-2)}$ (see the comment after Definition \ref{def:dmon}). This causes no extra difficulties, since $\psi^{(d-2)}(x)= -\int_x^{x^*} g(s)ds$ implies (\ref{eqn:nonstr_regv}) for $k=d-1$ upon changing variables to $t:=(x^*-s)^{-1}$. We now proceed by induction over $k$. The induction assumption is that, as $x\to\infty$, 
$$
\psi^{(k-1)}(x^*-x^{-1}) \sim c_{k-1}x^{-\alpha+k-1}L(x), \qquad -\alpha+k-1<0 \mathrm{\ and\ } c_{k-1}\neq 0.
$$
Here $L$ is some slowly varying function. The Monotone Density Theorem and (\ref{eqn:nonstr_regv}) then imply that $-x^{-2}\psi^{(k)}(x^*-x^{-1})\sim (\alpha-k+1)c_{k-1}x^{-\alpha+k-2} L(x)$, or, equivalently,
\begin{equation} \label{eqn:nonstr_regv2}
\psi^{(k)}(x^*-x^{-1}) \sim c_{k}x^{-\alpha+k}L(x),
\end{equation}
where $c_k=(-\alpha+k-1)c_{k-1}\neq 0$.

To complete the induction step, we must show that $-\alpha+k<0$. So suppose first that $-\alpha+k> 0$. Then (\ref{eqn:nonstr_regv2}) implies that $(-1)^k\psi^{(k)}(x)\to +\infty$ as $x\to x^*$, which contradicts $d$-monotonicity. Hence $-\alpha+k\leq 0$ for $k=1,\dots,d-1$, which in particular gives the last statement in part $(i)$.

We still have to exclude $-\alpha+k=0$. However, it is enough to show this for $k=1,\dots,d-2$, since this still would allow us to reach (\ref{eqn:nonstr_regv2}) for $k=d-1$ and $c_{d-1}\neq 0$, but possibly with $-\alpha+d-1=0$. So suppose $-\alpha+k=0$ for any $k\in\{1,\dots,d-2\}$. Then $(-1)^k\psi^{(k)}(x^*-x^{-1})\in RV_0$, but at the same time $d$-monotonicity implies that $(-1)^k\psi^{(k)}$ is convex for $k=1,\ldots,d-2$, and thus in particular
\begin{eqnarray*}
(-1)^k\psi^{(k)}\left(x^*-\frac{1}{2}x^{-1}\right) &=& (-1)^k\psi^{(k)}\left(\frac{1}{2}(x^*-x^{-1})+\frac{1}{2}x^*\right) \\
& \leq& (-1)^k\frac{1}{2}\psi^{(k)}(x^*-x^{-1})+(-1)^k\frac{1}{2}\psi^{(k)}(x^*) \\
&=& (-1)^k\frac{1}{2}\psi^{(k)}(x^*-x^{-1}),
\end{eqnarray*}
since $\psi^{(k)}(x^*)=0$. Equivalently,
\[ 2 \leq \frac{\psi^{(k)}(x^*-x^{-1})}{\psi^{(k)}(x^*-(2x)^{-1})}. \]
But here the right side tends to one as $x\rightarrow\infty$ since $(-1)^k\psi^{(k)}(x^*-x^{-1})\in RV_0$, so the inequality is eventually violated. This gives the desired contradiction.

The statement that was to be shown now follows upon noting that the induction assumption is satisfied for $k=1$ by hypothesis.

The proof of $(ii)$ is simpler. This time the induction assumption is 
$$
(-1)^k\psi^{(k)}(x^*-x^{-1})\in\RV_{-\alpha-(d-1)+k},
$$
where $k\leq d-1$. It clearly holds for $k=d-1$ by hypothesis. To prove the induction step, we use (\ref{eqn:nonstr_regv}) together with the induction assumption to obtain
$$
(-1)^{k-1}\psi^{(k-1)}(x^*-x^{-1}) = \int_x^\infty t^{-2}  (-1)^k \psi^{(k)}(x^*-t^{-1})dt = \int_x^\infty t^{-1-\alpha-d+k}L(t)dt
$$
for some slowly varying $L$. Since $\alpha+d-k>0$, Karamata's Theorem gives
$$
(-1)^{k-1}\psi^{(k-1)}(x^*-x^{-1}) \sim \frac{1}{\alpha+d-k} x^{-\alpha-(d-1)+(k-1)} L(x),
$$
and so $(-1)^{k-1}\psi^{(k-1)}(x^*-x^{-1})\in\RV_{-\alpha-(d-1)+(k-1)}$.
\proofend \end{proof}

Note that for $\alpha\in(0,\infty)$, Proposition \ref{prop:nonstrict_regvar} says in particular that $\psi(x^*-x^{-1})$ is in $\RV_{-\alpha-d+1}$ if and only if $\psi^{(d-1)}(x^*-x^{-1})$ is in $\RV_{-\alpha}$. We use this fact to prove the following connection between the generator $\psi$ and the radial part $R\sim F$. As in the case of rapid variation, Section \ref{section:rapvar}, we restrict ourselves to differentiable distribution functions.

\begin{proposition} \label{prop:nonstrict_regvar2}
Let $F$ be a distribution function with $F(0)=0$ and right endpoint $x^*<\infty$, let $\psi=\Wd F$ for $d\geq 2$ be its Williamson $d$-transform, and suppose that either $F$ or $\psi^{(d-1)}$ is differentiable. Then, for $\alpha>0$,
$$
\bar F(x^*-x^{-1}) \in \RV_{-\alpha} \quad \mathrm{if\ and\ only\ if} \quad \psi(x^*-x^{-1})\in \RV_{-\alpha+d-1}.
$$
\end{proposition}

\begin{proof}
We prove the ``only if'' part. By Lemma \ref{lemma:psiFprim}, both $F$ and $\psi^{(d-1)}$ are differentiable, and hence
$$
\bar F(x^*-x^{-1}) = -\int_x^\infty t^{-2}\bar F'(x^*-t^{-1})dt,
$$
since $\bar F(x^*)=0$. By the Monotone Density Theorem, $-\bar F'(x^*-x^{-1})\in\RV_{-\alpha+1}$. Lemma \ref{lemma:psiFprim} says that
$$
-\bar F'(x^*-x^{-1}) = \frac{1}{(d-1)!} (x^*-x^{-1})^{d-1}(-1)^d \psi^{(d)}(x^*-x^{-1}),
$$
which implies that also $(-1)^d\psi^{(d)}(x^*-x^{-1})$ is in $\RV_{-\alpha+1}$, since $(x^*-x^{-1})^{d-1}$ tends to a constant and thus is slowly varying. Now,
$$
(-1)^{d-1}\psi^{(d-1)}(x^*-x^{-1}) = \int_x^\infty t^{-2}(-1)^d\psi^{(d)}(x^*-t^{-1})dt,
$$
where the right side is equal to $\int_x^\infty t^{-\alpha-1}L(t)dt$ for some slowly varying $L$, due the regular variation of $(-1)^d\psi^{(d)}(x^*-x^{-1})$. Karamata's Theorem then implies that $(-1)^{d-1}\psi^{(d-1)}(x^*-x^{-1})\in\RV_{-\alpha}$, and we may apply Proposition \ref{prop:nonstrict_regvar} to reach the desired conclusion.

To prove the converse, one argues analogously, but in the reverse order. First one establishes regular variation of $\psi^{(d-1)}(x^*-x^{-1})$ using Proposition \ref{prop:nonstrict_regvar}. By the differentiability
assumption, one can differentiate once more and use the Monotone Density Theorem together with Lemma \ref{lemma:psiFprim} to obtain regular variation of $-\bar F'(x^*-x^{-1})$. Karamata's  theorem then gives the conclusion about $\bar F$.
\proofend \end{proof}

The intuition behind this result is that when $F$ (and hence $\psi$) has a finite right endpoint $x^*$, then the asymptotic behavior of $F$ near $x^*$ is the same as the asymptotic behavior of $\psi^{(d-1)}$ near $x^*$. This is different from the case when $x^*$ is infinite---then $F$ and $\psi$ have the same asymptotics, according to Proposition \ref{prop:Wd_tail_s1}.

\section{Behavior at the origin} \label{section:Wd_orig}

We conclude the investigation of the asymptotics of the Williamson $d$-transform by considering the properties near the origin.

\begin{proposition} \label{prop:3}
Let $R$ be a random variable on $(0,\infty)$ with Williamson $d$-transform $\psi$. Then
\begin{enumerate}
\item[(i)] If $1/R\in \RV_{-\alpha}$, $\alpha\in(0,1)$, then $1-\psi(x^{-1})\in\RV_{-\alpha}$.
\item[(ii)] If $1/R\in \RV_{-\alpha}$, $\alpha\in[1,\infty]$, then $1-\psi(x^{-1})\in \RV_{-1}$.
\item[(iii)] If $1-\psi(x^{-1})\in \RV_{-\alpha}$, $\alpha\in(0,1)$, then $1/R \in \RV_{-\alpha}$.
\end{enumerate}
\end{proposition}

\begin{proof}
Let $S_1$ be a ${\rm Beta}(1,d-1)$ distributed random variable independent of $R$. Since $S_1$ is continuous we have
$$
\P[1/S_1 > x] = \P[S_1 \leq   1/x] = 1 - (1-x^{-1})^{d-1} \sim x^{-1}, \qquad x\to\infty,
$$
which implies that $1/S_1\in\RV_{-1}$. Furthermore, it follows from the discussion after Proposition \ref{prop:Wd_prop} that $1-\psi(x^{-1})$ is the survival function of $1/(RS_1)$.

To show $(i)$, we pick an $\varepsilon>0$ such that $\alpha + \varepsilon < 1$, which is possible since $\alpha<1$. Then because $1/S_1\in\RV_{-1}$, $\E[(1/S_1)^{\alpha+\varepsilon}] <\infty$. Now, $1/R\in\RV_{-\alpha}$ by assumption, so Breiman's result, Lemma \ref{lemma:Breiman}, implies that $1/(RS_1)\in \RV_{-\alpha}$, or, in other words, $1-\psi(x^{-1})\in\RV_{-\alpha}$.

As for $(ii)$, consider first the case when $\alpha\in (1,\infty]$. Then because $1/R\in\RV_{-\alpha}$, $\E[(1/R)^{1+\varepsilon}] < \infty$ whenever $\varepsilon$ is such that $1+\varepsilon < \alpha$. One may then use Breiman's result again to conclude that $1/(RS_1)\in \RV_{-1}$, i.e. $1-\psi(x^{-1})\in\RV_{-1}$. 
If $\alpha=1$, Breiman's result is not applicable. Instead we show explicitly that in this case,
$$
L(x) := \frac{1-\psi(1/x)}{1/x}, \quad x\in(0,\infty),
$$
is slowly varying. To do so, first use the definition of the Williamson $d$-transform and the binomial expansion theorem to write 
$$
L(x) =  -\sum_{k=1}^{d-1} \binom{d-1}{k} (-1)^k x^{-k+1} \E[R^{-k} \I_{\{R >1/x\}}], \quad x\in (0,\infty).
$$
Further, denote by $\bar G$ the survival function of $1/R$ and  integrate by parts to obtain
$$
\E[R^{-k} \I_{\{R >1/x\}}] = - \int_0^x y^k d\bar G(y) = -x^k\bar G(x) + k\int_0^x \bar G(y) y^{k-1} dy.
$$
Thus $L(x)$ may be expressed as
$
L(x) = A(x) \{B_1(x) + \dots + B_{d-1}(x)\}
$
where
$$
A(x)  = \int_0^x \bar G(y) dy, \quad
B_1(x)  = (d-1)\left\{1- \frac{x\bar G(x)}{A(x)}\right\}
$$ 
and 
$$
B_k(x) = {d-1 \choose k} (-1)^{k-1} \frac{x \bar G(x)}{A(x)}\left(
-1 + \frac{k\int_0^x \bar G(y) y^{k-1} dy }{x^k \bar G(x)}\right), \quad k = 2,\dots, (d-1).
$$
Karamata's Theorem now implies that 
$$
\lim_{x\to\infty} \frac{x\bar G(x)}{A(x)} = 0,\quad \lim_{x\to\infty}\frac{k\int_0^x \bar G(y) y^{k-1} dy }{x^k \bar G(x)} = \frac{k}{k-1}, \quad k=2,\dots, (d-1).
$$
Thus as $x\to \infty$,  $B_1(x) \to d-1$ while $B_k(x) \to 0$ for $k=2,\dots, (d-1)$. Since $A(x)$ is slowly varying by Karamata's Theorem, $L(x)$ is a product of two slowly varying functions and therefore itself slowly varying.

Part $(iii)$ is based on the Monotone Density Theorem. Suppose $1-\psi(x^{-1}) = x^{-\alpha}L(x)$, with $L$ slowly varying. Since $\psi(0)=1$ we get
$$
1-\psi(x^{-1}) = -\int_0^{1/x}\psi'(y)dy  = -\int_x^\infty y^{-2} \psi'(y^{-1})dy.
$$
The Monotone Density Theorem gives $-x^{-1}\psi'(x^{-1})\sim\alpha x^{-\alpha}L(x)$ as $x\rightarrow\infty$. Suppose now that $(-1)^kx^{-k}\psi^{(k)}(x^{-1})\sim c_k x^{-\alpha} L(x)$ for some $1\leq k\leq d-2$, where $c_k$ is a constant. Since $\psi^{(k)}(x) = -\int_x^\infty \psi^{(k+1)}(y)dy$ we get, after changing variables,
$$
(-1)^k \psi^{(k)} (x^{-1}) = \int_0^x (-1)^{k+1}\psi^{(k+1)}(y^{-1})y^{-2}dy.
$$
By the Monotone Density Theorem, $(-1)^{k+1}\psi^{(k+1)}(x^{-1})x^{-2}\sim c_k (-\alpha+k)x^{-\alpha+k-1}L(x)$, and so by induction we obtain, for $k=0,1,\ldots,d-1$,
$$
(-1)^kx^{-k}\psi^{(k)}(x^{-1})\sim c_k x^\alpha L(x), \qquad x \rightarrow \infty,
$$
where the $c_k$ satisfy $c_0=1$, $c_1=\alpha$ and $c_{k+1}=(k-\alpha)c_k>0$. The inversion formula for the Williamson $d$-transform now yields
\begin{eqnarray*}
F(x^{-1})  &=& 1 - \psi(x^{-1}) - \sum_{k=1}^{d-1}\frac{(-1)^kx^{-k}\psi^{(k)}(x^{-1})}{k!}  \\
&\sim& \sum_{k=0}^{d-1} \frac{1}{k!} c_k x^{\alpha}L(x) = B (1-\psi(x^{-1})), \quad x\rightarrow\infty,
\end{eqnarray*}
where the constant $B$ is positive since each $c_k$ is.
\proofend \end{proof}

Suppose that $\psi$ is a $d$-monotone Archimedean generator. In particular, this implies that $\tilde \psi(x) := \psi(-x)$ is a distribution function on $(-\infty,0]$. Our final result shows that $d$-monotonicity of $\psi$ restricts the possible extremal behavior of $\tilde \psi$. In particular, it shows that Proposition \ref{prop:3} indeed covers (almost) all relevant cases (the exception being $\alpha=1$ in part $(iii)$, which is not covered).

\begin{proposition} \label{prop:2} Let  $\psi$ be a $d$-monotone Archimedean generator and $\tilde\psi$ as above. Then
\begin{enumerate}
\item[(i)] If $1-\psi(x^{-1})\in\RV_{-\alpha}$, then $\alpha\in[0,1]$. 
\item[(ii)] If the extremal behavior of $\tilde \psi$ is non-degenerate, then $\tilde \psi$ is in the maximum domain of attraction of the  Weibull distribution with the parameter range restricted to $(0,1]$. 
\end{enumerate}
\end{proposition}

\begin{proof}
Part $(i)$ stems from the fact that if $\alpha >1$, 
$$
\lim_{x\downarrow 0}\frac{1-\psi(x/2)}{1-\psi(x)} = 2^{-\alpha} < \frac{1}{2}. 
$$
In particular, therefore, there exists an $x_0 \in (0,\infty)$ such that $\psi(x/2) > \{\psi(0)+\psi(x)\}/2$ whenever $x\in(0,x_0)$. But $\psi$ is $d$-monotone and hence convex, so we obtain a contradiction.

Turning to $(ii)$, observe that since the upper endpoint of $\tilde\psi$ is finite, $\tilde\psi$ cannot be in the Fr\'echet domain of attraction. Furthermore, a necessary condition for $\tilde \psi\in\MDA(\Lambda)$ is that $1-\psi(x^{-1})\in\RV_{-\infty}$, which is excluded by $(i)$. So $\tilde \psi$ is in the Weibull domain.
This implies in particular that $1-\psi(x^{-1})\in\RV_{-\alpha}$ for some $\alpha\in(0,\infty)$. Due to $(i)$ however, $\alpha > 1$ is not possible.
\proofend \end{proof}

\section{Conditions on the inverse generator}

The most common approach found in the literature for dealing with Archimedean copulas is to express them in terms of the inverse generator $\phi = \psi^{-1}$. We have already seen some of the virtues of working with the generator $\psi$ instead, and the advantages will become even more apparent later on. However, in order to remain compatible with existing conventions we here give equivalent conditions on $\phi$ for the results in Sections \ref{section:regvar_strict}--\ref{section:Wd_orig}.

\begin{proposition} \label{prop:invgencond}
Let $\psi$ be an Archimedean generator, and let $\phi=\psi^{-1}$ be the inverse generator.
\begin{itemize}
\item[(i)] $\psi\in \RV_{-\alpha}$, $\alpha\in(0,\infty)$ if and only if $\phi(1/s)\in \RV_{1/\alpha}$, or, equivalently,
$$
\lim_{s\downarrow 0} \frac{\phi(cs)}{\phi(s)} = c^{-1/\alpha} \qquad \mathrm{for\ all\ } c>0.
$$
\item[(ii)] $\psi(x^*-x^{-1})\in \RV_{-\alpha}$, $\alpha\in [d-1,\infty)$, if and only if $(x^*-\phi(1/s))^{-1}\in\RV_{1/\alpha}$, or, equivalently,
$$
\lim_{s\downarrow 0} \frac{x^*-\phi(cs)}{x^*-\phi(s)} = c^{1/\alpha}\qquad  \mathrm{for\ all\ } c>0.
$$
Here $x^*$ is the right endpoint of $\psi$.
\item[(iii)] $1-\psi(x^{-1})\in\RV_{-\alpha}$, $\alpha\in(0,1]$, if and only if $\phi(1-s^{-1})\in\RV_{-1/\alpha}$, or, equivalently,
$$
\lim_{s\downarrow 0} \frac{\phi(1-cs)}{\phi(1-s)} = c^{1/\alpha}\qquad  \mathrm{for\ all\ } c>0.
$$
\item[(iv)] $\psi\in\MDA(\Lambda)$ if and only if $\phi(1/s)$ is in de Haan's $\Pi$-class\footnote{A non-negative, non-decreasing function $g$ is in de Haan's $\Pi$-class if $\lim_{s\to\infty}\{g(ts)-b(s)\}/c(s) = \log t$, $t\in\R$, for some functions $b(s)\in\R$ and $c(s)>0$. The function $c(s)$ is called the \emph{auxiliary function} of $g$. See Resnick \cite{Resnick:2008}, Section 0.4.3 for an concise account of the relevant theory.}. The auxiliary function $a(x)$ of $\psi(x)$ and $c(s)$ of $\phi(1/s)$ are related through $c(s)=a(\phi(1/s))$ and $a(x) = c(1/\psi(x))$.
\end{itemize}
\end{proposition}

To simplify matters, we need the following adaptation of a theorem on inverses of regularly varying functions:

\begin{lemma} \label{lemma:aux1}
Suppose $f$ is invertible and $\alpha\in(0,\infty)$. Then $f\in \RV_{-\alpha}$, $\alpha>0$ if and only if $f^{-1}(s^{-1})\in\RV_{1/\alpha}$, or, equivalently,
$$
\lim_{s\downarrow 0} \frac{f^{-1}(cs)}{f^{-1}(s)} = c^{-1/\alpha}, \qquad \mathrm{for\ all\ } c >0.
$$
\end{lemma}

\begin{proof}
$f\in\RV_{-\alpha}$ if and only if $1/f\in\RV_{\alpha}$. An easy computation shows that $[1/f]^{-1}(s) = f^{-1}(s^{-1})$, which by Theorem 1.5.12 in Bingham et al.~\cite{Bingham/Goldie/Teugels:1987} is equivalent to $f^{-1}(s^{-1})\in\RV_{1/\alpha}$. The last equivalent formulation follows directly from the definition of regular variation upon changing the variable from $s$ to $s^{-1}$.
\proofend \end{proof}

\begin{proof} \emph{(of Proposition \ref{prop:invgencond})}
Part $(i)$ follows directly from Lemma \ref{lemma:aux1} since $\psi$ is strict and hence invertible. For Part $(ii)$, note that $x^*-x^{-1}\in\{\psi>0\}$. On this set $\psi$ is invertible, and the inverse of $\psi(x^*-x^{-1})$ is $1/(x^*-\phi(s^{-1}))$. Now apply Lemma \ref{lemma:aux1}. Part $(iii)$ follows in the same fashion since $\psi(x^{-1})$ is invertible for $x$ sufficiently large, with inverse $\phi(1-s^{-1})$. For Part $(iv)$, Proposition 0.10 in Resnick \cite{Resnick:2008} says that $\psi\in\MDA(\Lambda)$ if and only if $(1/\psi)^{\leftarrow}(s) = \inf\{x\geq 0 \mid 1/\psi(x)\geq s\}$ is in de Haan's $\Pi$-class. But since $\psi$ is continuous, we have $\inf\{x\geq 0 \mid 1/\psi(x)\geq s\} = \inf\{x\geq 0 \mid 1/\psi(x) = s \}$, which is equal to $\phi(s^{-1})$ by definition. The relationship between the auxiliary functions follows from Proposition 0.9 in \cite{Resnick:2008}, where, in order to obtain the expression for $a(x)$, we used that $[\phi(1/s)]^{\leftarrow} = 1/\psi$ since $\phi$ is a right inverse of $\psi$: $\phi(1/s)=x$ implies $1/s=\psi(x)$.
\proofend \end{proof}

\chapter{Coefficients of tail dependence}  \label{section:taildep}

\section{Coefficient of the first kind}

An important problem in probability theory and statistics is to find simple measures, preferably scalar, that quantify the strength of dependence between random variables. In the context of extremes, one may ask about the likelihood of observing large values in one component $X_1$ given that another component $X_2$ is large. One such measure is the \emph{coefficient of (upper) tail dependence $\lambda_u$}, defined as
$$
\lambda_u := \lim_{q\to 1} \P[X_1 > F_1^{\leftarrow}(q) \mid X_2 > F_2^{\leftarrow}(q)],
$$
where $F_1$ and $F_2$ are the marginal distribution functions of $X_1$ and $X_2$, respectively, and the limit is assumed to exist (this is not always the case!). The \emph{coefficient of lower tail dependence}, $\lambda_l$, is defined analogously:
$$
\lambda_l := \lim_{q\to 0} \P[X_1 \leq F_1^{\leftarrow}(q) \mid X_2 \leq F_2^{\leftarrow}(q)].
$$
Note that dependence is measured on a quantile scale. It is therefore not surprising that these measures can be expressed in terms of the associated copula $C$. Indeed, using the definition of the generalized inverse, one easily checks the validity of
$$
 \P[X_1 \leq F_1^{\leftarrow}(q) \mid X_2 \leq F_2^{\leftarrow}(q)] = \frac{\P[F_1(X_1)\leq q, F_2(X_2)\leq q]}{\P[F_2(X_2)\leq q]} = \frac{C(q,q)}{q}.
$$
Similarly one can show that $\lambda_u$ also depends only on the copula. In the case where $C$ is Archimedean, $\lambda_u$ and $\lambda_l$ have simple forms:

\begin{proposition} \label{prop:taildep_psi}
Provided the limits exist, the coefficients of upper and lower tail dependence of an Archimedean copula with generator $\psi$ are given by
\begin{eqnarray*}
\lambda_u &=& 2-\lim_{t\rightarrow 0} \frac{1-\psi(2t)}{1-\psi(t)} \\
\lambda_l &=&\left\{ \begin{array}{ll}
  \lim_{t\rightarrow\infty}  \frac{\psi(2t)}{\psi(t)}, & \mathrm{if \ }\psi \ \mathrm{is\ strict}   \\
  0,&  \mathrm{otherwise}
\end{array}\right.
\end{eqnarray*}
\end{proposition}

\begin{proof}
Let $(U_1,\dots,U_d)$ be distributed according to $C$. Due to exchangeability, all two-dimensional marginals are given by $C(u_1,u_2,1,\ldots,1)=\psi(\psi^{-1}(u_1)+\psi^{-1}(u_2))$. The expression for $\lambda_u$ and the fact that the marginals are uniform yield
\begin{eqnarray*}
\lambda_u &=& \lim_{q\rightarrow1} \frac{ \P[U_1>q, U_2>q]}{\P[U_2>q]} \\
&=& \lim_{q\rightarrow1} \frac{1-\P[U_1\leq q]-\P[U_2\leq q]+\P[U_1\leq q,U_2\leq q]}{1-\P[U_2\leq q]} \\
&=& 2 - \lim_{t\rightarrow0}\frac{1-\psi(2t)}{1-\psi(t)}.
\end{eqnarray*}
Here we used that $q=\psi \circ \psi^{-1}(q)$ and that $t:=\psi^{-1}(q)\rightarrow 0$ as $q\rightarrow1$. Similarly,
\begin{eqnarray*}
\lambda_l &=&  \lim_{q\rightarrow 0} \frac{ \P[U_1\leq q, U_2\leq q]}{\P[U_2\leq q]} = \lim_{q\rightarrow0} \frac{\psi(2\psi^{-1}(q))}{q}= \lim_{t\rightarrow\psi^{-1}(0)}\frac{\psi(2t)}{\psi(t)}.
\end{eqnarray*}
If $\psi$ is strict, then $\psi^{-1}(0)=\infty$ and we get the result. If $\psi$ is non-strict, then $c := \psi^{-1}(0) < \infty$. So for all $t\in [c/2,c)$, the numerator $\psi(2t)$ is zero, while the denominator is non-zero. Hence the limit must be zero.
\proofend \end{proof}

Note that the generator in Example \ref{ex:no_MDA} provides an example where the index of lower tail dependence does not exist. However, the results from Section \ref{section:Wd_tail} allow us to give sufficient conditions for the existence of the coefficients of tail dependence and compute their values in terms of the corresponding radial part $R$.

\begin{corollary} \label{cor:taildep_psi}
Suppose the Archimedean copula $C$ has generator $\psi$, and let $R\sim F=\Wd^{-1}\psi$ be the corresponding radial distribution.
\begin{itemize}
\item[(i)] If $1/R\in\RV_{-\alpha}$, $\alpha\in(0,\infty)$, then $\lambda_u=2-2^{\beta}$, where $\beta=\alpha \wedge 1$.
\item[(ii)] If $R\in\RV_{-\alpha}$, $\alpha\in[0,\infty)$, then $\lambda_l=2^{-\alpha}$.
\end{itemize}
\end{corollary}

\begin{proof}
For $(i)$, we observe that $1-\psi(x^{-1})\in\RV_{-\beta}$ according to Proposition \ref{prop:3}, and hence
$$
\lim_{t\to 0}\frac{1-\psi(2t)}{1-\psi(t)} = \lim_{x\to \infty}\frac{1-\psi((x/2)^{-1})}{1-\psi(x^{-1})} = 2^{\beta},
$$
by definition of regular variation. The expression in Proposition \ref{prop:taildep_psi} for $\lambda_u$ gives the result. For $(ii)$, the result is immediate from Proposition \ref{prop:Wd_tail_s1}, the definition of regular variation, and Proposition \ref{prop:taildep_psi}.
\proofend \end{proof}

We remark that this result does not give conditions on $R$ ensuring that $\lambda_u=1$. This is due to the limitation $\alpha>0$ in Proposition \ref{prop:3} $(i)$.

If $\psi\in\RV_{-\infty}$, then Proposition \ref{prop:taildep_psi} clearly implies that $\lambda_l=0$. Proposition \ref{prop:Wd_tail_Gumbel} gives one sufficient condition for this to happen, namely that the corresponding radial part $R$ is in the Gumbel domain of attraction and has differentiable distribution function with eventually monotone derivative.

Note that the converse of Corollary \ref{cor:taildep_psi} does not hold---knowing e.g.~$\lambda_l$ only gives information about the limit of $\psi(cx)/\psi(x)$ for $c=2$, which is not enough to deduce regular variation of $\psi$ (and thus of $R$). The situation is analogous for $\lambda_u$.

A natural extension of the coefficients of tail dependence for high-dimensional distributions is
$$
\lambda^{m,n}_l := \lim_{q\to 0} \P[X_1 \leq F_1^{\leftarrow}(q), \dots, X_m \leq F_m^{\leftarrow}(q) \mid X_{m+1} \leq F_{m+1}^{\leftarrow}(q), \dots, X_{m+n} \leq F_{m+n}^{\leftarrow}(q)],
$$
and analogously for $\lambda^{m,n}_u$. As in Proposition \ref{prop:taildep_psi} one shows that for Archimedean copulas, we have the expression
$$
\lambda^{m,n}_l =\lim_{t\to\infty} \frac{\psi((m+n)t)}{\psi(nt)}
$$
if $\psi$ is strict. If $\psi$ is regularly varying with index $-\alpha$, $\alpha\in [0,\infty)$ (which happens if $R$ has the same property), we conclude as before that
\begin{equation} \label{eqn:tcmv}
\lambda^{m,n}_l = \left(\frac{m}{n}+1\right)^{-\alpha}.
\end{equation}
The expression for $\lambda_u$ is possible to derive, but has a more complicated form since its derivation requires that survival functions of the copula $C$ be expressed in terms of marginal distribution functions of different dimension. This can be done using the so-called sieve formula, but the calculations are cumbersome.

As a final observation we note that although there is no hope of a finding a converse to Corollary \ref{cor:taildep_psi}, the situation is in a sense different for $\lambda_l^{m,n}$. For, if we know that $\lambda_l^{m,n}$ has the form (\ref{eqn:tcmv}) for all $m,n\in \mathbb N$, then $\psi$ must be regularly varying with index $-\alpha$. This follows from Theorem 1.4.3 in Bingham et al. \cite{Bingham/Goldie/Teugels:1987}, which says that $\psi\in\RV_{-\alpha}$ if and only if $\psi(cx)/\psi(x)\to c^{-\alpha}$ for all $c$ in a dense subset of $(0,\infty)$, provided the condition
$$
\limsup_{c\downarrow 1} \psi^*(c) \leq 1, \quad \mathrm{where\ } \psi^*(c):=\limsup_{x\to\infty} \psi(cx)/\psi(x),
$$
is satisfied. But this automatically holds since $\psi$ is decreasing and hence $\psi^*(c)\leq 1$ for all $c\geq 1$. Furthermore, although (\ref{eqn:tcmv}) only gives $\psi(cx)/\psi(x)\to c^{-\alpha}$ for all rational $c$ in $(1,\infty)$, this is enough since for $0<c<1$ rational, it implies
$$
\lim_{x\to\infty} \psi(cx)/\psi(x) = \lim_{x\to\infty} \psi(x)/\psi(x/c) = \left\{ \lim_{x\to\infty} \psi(x/c)/\psi(x) \right\}^{-1} = c^{-\alpha}.
$$
Of course, this way of inferring properties of $\psi$, and hence $R$, based on the behavior of $\lambda_l^{m,n}$ is of little practical importance since \emph{all} these coefficients must satisfy (\ref{eqn:tcmv}). The corresponding estimation problem is obviously infeasible. Moreover, the interpretation of $\lambda_l^{m,n}$ as describing multivariate marginal tail behavior only remains valid if $\psi$ generates a copula in \emph{any} dimension, and thus is completely monotone.

\section{Coefficient of the second kind}

When a coefficient of tail dependence is zero, this is interpreted as \emph{asymptotic independence}\footnote{This term has several different meanings depending on the context. Here we use it a vague sense to support the intuition about the coefficients of tail dependence.} of the tails of the marginals, in the sense that unlikely events tend to occur more or less independently of each other. The ``more or less'' part is of course not completely satisfactory, and supplementary measures of tail dependence have been suggested. One such measure was introduced by Coles et al.~\cite{Coles/Heffernan/Tawn:2000}, and following Demarta \cite{Demarta:2007} we call it the \emph{(upper) coefficient of tail dependence of the second kind}. It is defined as
$$
\bar\lambda_u = \lim_{q\to 0}\frac{\log \left( \P[X_1 > F_1^{\leftarrow}(q)] \P[X_2 > F_2^{\leftarrow}(q)] \right)}{\log\left( \P[X_1 > F_1^{\leftarrow}(q), X_2 > F_2^{\leftarrow}(q)] \right)} - 1,
$$
with the lower coefficient $\bar\lambda_l$ being defined analogously. These quantities compare, on a logarithmic scale, the joint exceedance probabilities under the actual dependence structure with the exceedance probabilities one would observe if the components were independent. Expressed in terms of the copula $C$ of $(X_1,X_2)$, the tail dependence coefficients of the second kind are
\begin{equation} \label{eqn:taildep2k}
\bar\lambda_u = \lim_{q\to 1} \frac{2\log (1-q)}{\log (1-2q+C(q,q))} - 1 \qquad \mathrm{and} \qquad \bar\lambda_l = \lim_{q\to 0} \frac{2\log q}{\log C(q,q)} - 1.
\end{equation}
We first consider the lower coefficient $\bar \lambda_l$. It turns out that when $\psi\in\RV_{-\alpha}$ for $\alpha>0$, and hence $\lambda_l>0$, $\bar \lambda_l$ gives no additional information:

\begin{proposition} \label{prop:123}
Suppose $\psi\in\RV_{-\alpha}$, $\alpha>0$. Then $\bar\lambda_l=1$.
\end{proposition}

\begin{proof}
Change the variable to $x:=\psi^{-1}(q)$ in (\ref{eqn:taildep2k}) and apply L'H\^opital's rule:
$$
 \lim_{q\to 0} \frac{2\log q}{\log C(q,q)} =  \lim_{x\to \infty} \frac{2\log \psi(x)}{\log \psi(2x)} = \lim_{x\to \infty} \frac{\psi(2x)}{\psi(x)} \frac{\psi'(2x)}{\psi'(x)}.
$$
Now, $\psi'$ is monotone by $d$-monotonicity, and $\psi(x)= - \int_x^\infty \psi'(t)dt$. Since $\alpha\neq 0$, the Monotone Density Theorem implies that $-\psi' \in \RV_{-\alpha-1}$, and so the limit above is equal to $2^{-\alpha}/2^{-\alpha-1}=2$. By (\ref{eqn:taildep2k}), then, $\bar \lambda_l = 1$.
\proofend \end{proof}

In fact one can easily show that $\bar \lambda_l$ is always one whenever $\lambda_l>0$ (see Demarta \cite{Demarta:2007}, Lemma 2.1 $(iii)$ for the corresponding property in the case of $\bar \lambda_u$). If $C$ is Archimedean with strict generator $\psi$, and $\psi$ is in the Gumbel domain of attraction (and hence in $\RV_{-\infty}$, implying that $\lambda_l=0$), we can be more precise about the behavior of $\bar\lambda_l$.

\begin{proposition} \label{prop:barla}
Suppose $C$ is a bivariate Archimedean copula whose generator $\psi$ is strict and in the Gumbel domain of attraction with auxiliary function $a$. If $a\in\RV_{\beta}$, then
$$
\bar\lambda_l = 2^{-\beta}-1.
$$
\end{proposition}

\begin{proof}
As in the preceding proof we obtain
$$
\lim_{q\to 0} \frac{2\log q}{\log C(q,q)} =  \lim_{x\to\infty} \frac{\psi'(x)/\psi(x)}{\psi'(2x)/\psi(2x)}.
$$
By Theorem \ref{thm:MDA_char} and the subsequent remark, we have $a(x)\sim -\psi(x)/\psi'(x)$, whence the result follows due to the assumption on $a$ and the definition of regular variation.
\proofend \end{proof}

These results are in line with the intention that $\bar\lambda_l$ should allow to differentiate between different degrees of asymptotic independence (in the above vague sense). The conditions on $\psi$ in Proposition \ref{prop:barla} can immediately be translated into conditions on the radial part $R$ using Proposition \ref{prop:Wd_tail_Gumbel}.

As mentioned above, there is a counterpart to Proposition \ref{prop:123} for the upper coefficient $\bar \lambda_u$. More precisely, $\bar \lambda_u = 1$ as soon as $\lambda_u >0$. A sufficient condition for this to happen is, according to Corollary \ref{cor:taildep_psi}, that the radial part $R$ corresponding to $\psi$ satisfies $1/R\in\RV_{-\alpha}$ with $\alpha\in(0,1)$. We have not been able to find a counterpart to Proposition \ref{prop:barla}, i.e.~to deduce what happens in the case $\alpha\in[1,\infty)$, and hence $\lambda_u=0$.

We conclude this chapter with a short discussion on the relationship between distributions, their copulas, survival copulas, and coefficients of tail dependence. The tail dependence coefficients are sometimes defined as a quantity relating to the \emph{distribution} $F$ of a random vector, but having the property that it only depends on the (ordinary or survival) copula of that distribution. That is how the concept was approached here. However, many authors define the tail dependence coefficients as quantities relating directly to a particular \emph{copula} $C$, regardless of whether this copula in a specific modeling context appears as the (ordinary) copula of a distribution under consideration, or as its survival copula. This latter convention has the consequence that if $\lambda_u$, say, is the upper tail dependence coefficient of $C$, and $C$ is the (ordinary) copula of $F$, then $\lambda_u$ describes the \emph{upper} tail of $F$. If, on the other hand, $C$ is the survival copula of $F$, then $\lambda_u$ describes the \emph{lower} tail of $F$. When dealing with Archimedean copulas, this can become confusing since these copulas appear as the \emph{survival copulas} of $\ell_1$-norm symmetric distributions. Hence the quantities $\lambda_l$ and $\bar\lambda_l$ treated in this chapter both describe the \emph{upper} tail of the corresponding $l_1$-norm symmetric distribution, while $\lambda_u$ and $\bar\lambda_u$ describe its \emph{lower} tail.

\chapter{Extreme value limits of Archimedean copulas} \label{section:EV_limits}

A $d$-dimensional copula $C$ with the property that $C^n(\bu^{1/n})=C(\bu)$ for all $n>0$ and $\bu\in[0,1]^d$  is called \emph{extreme value} (EV) copula. These copulas arise when one considers maxima (or, as we will see, minima) of multivariate distributions. It is well-known (see Genest and Rivest \cite{Genest/Rivest:1989}) that if $C$ is Archimedean, it belongs to the Gumbel family, i.e.~its generator is $\psi^{\mathrm{Gu}}_\theta(x)=\exp(-x^{1/\theta})$ for some $\theta\geq 1$. The language of $\ell_1$-norm symmetric distributions makes it simple to understand why this is to be expected. This will be apparent in Section \ref{section:L1nsd_min}, where we consider the limiting distributions of componentwise minima of $\ell_1$-norm symmetric distributions. In Section \ref{section:L1nsd_max} we consider the maxima of these distributions. Here the limit turns out to be the so-called Galambos copula (see \cite{Galambos:1975} and \cite{Caperaa/Fougeres/Genest:2000}).

\section{Componentwise minima of $\ell_1$-norm symmetric distributions} \label{section:L1nsd_min}

Consider an i.i.d.~sequence $\bX_1,\bX_2,\ldots$ of $\R^d$-valued random vectors. The corresponding \emph{componentwise $n$-block minima} are defined as
$$
\bW_n = (W_{1,n},\ldots,W_{d,n}) := \left(\min\{X_{1,k}\}_{k=1}^n,\ldots,\min\{X_{d,k}\}_{k=1}^n\right).
$$

In the sequel it will always be assumed that each $\bX_n=(X_{1,n},\ldots,X_{d,n})$ follows a $d$-dimensional $\ell_1$-norm symmetric distribution with radial part $R \sim F$, where $R>0$ a.s., and joint survival survival function $\bar H$. By Theorem \ref{thm:main} it then follows that $\bar H(\bx) = \psi(\|\bx\|_1)$ for $\bx\in\R^d_+$, where $\psi=\Wd F$. Moreover, $\psi$ is then the generator of the survival copula $C$ of $\bX_n$, which is Archimedean. We now fix this notation and do not repeat it at each instance it appears in the sequel.

The $\ell_1$-norm symmetric distributions have the following closure property under the operation of forming minima:

\begin{lemma} \label{lemma:min_L1nsd}
The $n$-block minima $\bW_n$ are again $\ell_1$-norm symmetric. The corresponding survival copula generators are $\psi^n$.
\end{lemma}

\begin{proof} 
For each $n$, $\boldsymbol W_n$ is clearly confined to the interior of $\R^d_+$ with probability one. Its survival function is
\begin{eqnarray*}
\P[\bW_n>\bx] &=&\P[\{\min\{X_{1,k}\}_{k=1}^n>x_1\}\cap\cdots\cap\{\min\{X_{d,k}\}_{k=1}^n>x_d\}] \\
&=& \P\left[ \bigcap_{i=1,\ldots,d} \bigcap_{k=1,\ldots,n} \{ X_{i,k}>x_i \} \right] \\
&=& \left( \P\left[ \bigcap_{i=1,\ldots,d} \{ X_{i,1}>x_i \} \right]  \right)^n  \\
&=& \psi^n(\|\bx\|_1),
\end{eqnarray*}
where the third equality uses mutual independence of $\bX_1,\ldots,\bX_n$. The result now follows by Theorem \ref{thm:main}.
\proofend \end{proof}

\begin{lemma} \label{lemma:1dim_min_conv}
Suppose $1-\psi(x^{-1})\in RV_{-1/\theta}$ for some $\theta<\infty$. Then there exists a sequence $(c_n)$ of positive constants such that $c_n^{-1}W_{1,n} \Rightarrow Y$ as $n\rightarrow\infty$, where the survival function of $Y$ is $\psi^{\mathrm{Gu}}_\theta$. 
\end{lemma}

\begin{proof} 
Using the identity $\max(-x)=-\min(x)$ we get $-W_{1,n}=\max\{-X_{1,k}\}_{k=1}^n$. The random variables $-X_{1,k}$ are i.i.d.~with distribution function $\psi(-x)$ for $x\leq 0$, which has right endpoint zero. Moreover, since $1-\psi(x^{-1})\in RV_{-1/\theta}$ by assumption,  $\psi(-x)$ is in the Weibull domain of attraction by Theorem \ref{thm:MDA_char}. By Proposition \ref{prop:normseq}, there exist positive constants $c_n$ such that $-c^{-1}_n W_{1,n} \Rightarrow -Y \sim \exp(-(-x)^{1/\theta})=\psi^{\mathrm{Gu}}_\theta(-x)$. The result now follows upon noting that, for $x\geq0$,
$$
\P[Y>x] = \P[-Y\leq -x] = \psi^{\mathrm{Gu}}_\theta(-(-x)),
$$
where we used the continuity of $\psi^{\mathrm{Gu}}_\theta$.
\proofend \end{proof}

Before continuing we note that Proposition \ref{prop:2} implies that $1/\theta\in[0,1]$. Since the hypothesis of Lemma \ref{lemma:1dim_min_conv} excludes $\theta=\infty$, we always have $\theta\in[1,\infty)$, which is consistent with the requirements of the Gumbel generator.

We are now in a position to prove the next result on the limiting distribution of normalized $n$-block minima of $\ell_1$-norm symmetric distributions. As a corollary we obtain the extreme value limits of Archimedean copulas.

\begin{proposition} \label{prop:min_limit}
Suppose $1-\psi(x^{-1})\in RV_{-1/\theta}$ for some $\theta<\infty$. Then there exists a sequence $(c_n)$ of positive constants such that $c_n^{-1}\bW_n \Rightarrow \bW_0$ as $n\rightarrow\infty$, where $\bW_0$ is $\ell_1$-norm symmetric with survival copula generator $\psi^{\mathrm{Gu}}_\theta$.
\end{proposition}

\begin{proof}
For any constant $a>0$, Lemma \ref{lemma:min_L1nsd} implies that the survival function of $a^{-1}\bW_n$ is $\psi^n(a \|\bx\|_1)$, with the one-dimensional marginals thus being $\psi^n(a x)$. With $c_n$ as in Lemma \ref{lemma:1dim_min_conv}, we therefore conclude that the functions $x\mapsto \psi^n(c_n x)$ converge weakly to $\psi^{\mathrm{Gu}}_\theta$ as $n\rightarrow\infty$. Hence $\P[c_n^{-1} \bW_n > \bx] = \psi^n(c_n \|\bx\|_1) \rightarrow \psi^{\mathrm{Gu}}_\theta(\|\bx\|_1)$ weakly as $n\rightarrow\infty$, which means precisely that $c_n^{-1}\bW_n \Rightarrow \bW_0$. Theorem \ref{thm:main} implies that $\boldsymbol W_0$ has the desired properties.
\proofend \end{proof}

\begin{corollary} \label{cor:Arch_limit}
Suppose $C$ is an Archimedean copula whose generator $\psi$ satisfies $1-\psi(x^{-1})\in RV_{-1/\theta}$ for some $\theta<\infty$. Then $C^n(\bu^{1/n})\rightarrow C^{Gu}_\theta(\boldsymbol u)$ pointwise.
\end{corollary}

\begin{proof}
Let the $\bX_n$ above follow the $\ell_1$-norm symmetric distribution whose survival copula generator is $\psi$. By Lemma \ref{lemma:min_L1nsd}, the survival copula of $\bW_n$, and hence of $a\bW_n$ for any $a>0$, is
\[ \psi^n(\psi^{-1}(u_1^{1/n})+\ldots+\psi^{-1}(u_1^{1/n}))=C^n(\bu^{1/n}). \]
Weak convergence of a sequence of joint distributions toward a limiting distribution implies pointwise convergence of the survival copulas toward the survival copula of the limiting distribution (see e.g.~Feidt et al.~\cite{Feidt/Genest/Neslehova:2008}). By Proposition \ref{prop:min_limit}, the limiting distribution has survival copula $C^{Gu}_\theta$, with $\theta$ as in the statement of the corollary. The result follows. 
\proofend \end{proof}

In the case of $\ell_1$-norm symmetric distributions, it is the \emph{survival copula} that is Archimedean. This is the reason why we considered the \emph{minima} of these distributions in order to arrive at the conclusion in  Corollary \ref{cor:Arch_limit}. However, in a situation where an Archimedean copula $C$ is the \emph{copula} (as opposed to the survival copula) of a $d$-dimensional distribution $\boldsymbol Y$ under consideration, then $C^n(\boldsymbol u^{1/n})$ will be the copulas of the \emph{$n$-block maxima}, defined as
$$
\boldsymbol{M}_n = (M_{1,n},\ldots,M_{d,n}) := \left(\max\{Y_{1,k}\}_{k=1}^n,\ldots,\max\{Y_{d,k}\}_{k=1}^n\right),
$$
where $\boldsymbol Y_n = (Y_{1,n},\dots,Y_{d,n})$, $n\geq 1$, form a sequence of independent copies of $\boldsymbol Y$. Indeed, if $G$ is the joint distribution function of $\boldsymbol Y$, independence of the $\boldsymbol M_n$ implies (cf.~the proof of Lemma \ref{lemma:min_L1nsd}) that
$$
\P[\boldsymbol M_n \leq \bx] = \left( \P[\boldsymbol Y \leq \bx] \right)^n = G^n(\bx).
$$
Hence the marginals are $\P[M_{i,n}\leq x_i]=G_i^n(x_i)$, and an explicit computation yields the form $C^n(\boldsymbol u^{1/n})$ of the copula of $\boldsymbol M_n$. Corollary \ref{cor:Arch_limit} of course applies also in this case. To summarize, we have the following (as before, $C$ is Archimedean with generator $\psi$):
\begin{itemize}
\item[(i)] If $C$ is the \emph{copula} of the distribution $\boldsymbol Y$, then the dependence structure of the corresponding $n$-block \emph{maxima} can be approximated by the Gumbel copula for large $n$ (if $1-\psi(x^{-1})\in \RV_{-\alpha}$).
\item[(ii)] If $C$ is the \emph{survival copula} of the distribution $\boldsymbol Y$, then the dependence structure of the corresponding $n$-block \emph{minima} can be approximated by the Gumbel copula for large $n$ (if $1-\psi(x^{-1})\in \RV_{-\alpha}$).
\end{itemize}
We finally remark that Proposition \ref{prop:3} allows us to express the condition on $\psi$ in terms of a condition on the corresponding radial part $R$---namely, $1/R$ must have a regularly varying tail.

\section{Componentwise maxima of $\ell_1$-norm symmetric distributions} \label{section:L1nsd_max}

In this section, $\boldsymbol M_n$ will denote the componentwise $n$-block maxima of the i.i.d.~$\ell_1$-norm symmetric random vectors $\bX_1,\bX_2,\ldots$, i.e.
$$
\boldsymbol{M}_n = (M_{1,n},\ldots,M_{d,n}) := \left(\max\{X_{1,k}\}_{k=1}^n,\ldots,\max\{X_{d,k}\}_{k=1}^n\right).
$$
As before, $R\sim F$ with $R>0$ a.s.~is the radial part of $\bX_n$, and $\psi=\Wd F$ is the corresponding Archimedean generator. The generated survival copula is denoted $C$. Unlike the minima $\boldsymbol{W}_n$ considered in Section \ref{section:L1nsd_min}, the $\boldsymbol{M}_n$ are \emph{not} $\ell_1$-norm symmetric. We thus need to take a different approach to analyze their asymptotic behavior, and it turns out that the theory of multivariate regular variation and multivariate extremes provides the proper framework. The stochastic representation of the $\ell_1$-norm symmetric distributions renders the computations particularly simple.

Since the multivariate extreme value theory is more complicated than its univariate counterpart, and because it is only used in this section, we refrain from giving an introduction to the general theory. Instead we give precise references for results that are used. The book \cite{Resnick:2008} by Resnick contains all theory we need here, and much more beyond that.

\begin{proposition} \label{prop:mv_regvar}
Suppose $\bX$ is $\ell_1$-norm symmetric with radial part $R\in\RV_{-\alpha}$, $\alpha\in(0,\infty)$. Then $\bX$ is in the domain of attraction of
$$
G(\bx) := \exp\left\{-\mu([\0,\bx]^c)\right\}, \qquad \bx>\0.
$$
where the \emph{exponent measure} $\mu$ is given by
$$
\mu\{\bx\in\R^d_+ \backslash \{\0\} \ | \ \|\bx\|_1>r, \bx/\|\bx\|_1\in A\} = r^{-\alpha}\sigma(A), \qquad r>0, A\in\mathcal{S}^d,
$$
and the \emph{spectral measure} $\sigma$ is the uniform distribution on $\mathcal{S}^d$.
\end{proposition}

\begin{proof}
Since $\bX$ is $\ell_1$-norm symmetric, Proposition \ref{prop:L1_properties} implies that $\boldsymbol{S}:\eqdis \bX/\|\bX\|_1$ is uniformly distributed on $\mathcal{S}^d$, $\|\bX\|_1\eqdis R$, and $\|\bX\|_1$ and $\bX/\|\bX\|_1$ are independent. Hence for any Borel subset $A\subseteq\mathcal{S}^d$ and $x>0$,
$$
\lim_{u\rightarrow\infty} \frac{\P[\|\bX\|_1>ux, \bX/\|\bX\|_1 \in A]}{\P[\|\bX\|_1>u]} = \lim_{u\rightarrow\infty} \frac{\bar{F}(ux)}{\bar{F}(u)}\P[\boldsymbol{S}\in A] = x^{-\alpha}\P[\boldsymbol{S}\in A].
$$
Setting $t:=1/\bar F(u)$ and $a_t:=(1/\bar F)^{\leftarrow}(t)$, we also have
$$
\lim_{u\to\infty} \frac{\P[\|\bX\|_1>ux, \bX/\|\bX\|_1 \in A]}{\P[\|\bX\|_1>u]} = \lim_{t \to\infty} t \P[a_t^{-1}\|\boldsymbol X\|_1>\bx, \ \bX/\|\bX\|_1\in A].
$$
We can now apply Corollary 5.18 (b) in Resnick \cite{Resnick:2008} to reach the desired conclusion with  $\sigma=\P[\boldsymbol{S}\in (\cdot)]$.
\proofend \end{proof}

This means in particular that $\bX$ is \emph{multivariate regularly varying} with tail index $\alpha$ and spectral measure $\sigma$. We now express the exponent measure $\mu$ directly in terms of $\sigma$.

\begin{lemma} \label{lemma:mu_repr}
Let $\mu$ and $\sigma$ be as in Proposition \ref{prop:mv_regvar}. Then 
$$
\mu([\0,\bx]^c) = \int_{\mathcal{S}^d} \max \left\{\frac{w^\alpha_1}{x_1^\alpha},\ldots,\frac{w^\alpha_d}{x_d^\alpha}\right\}\sigma(d\boldsymbol{w}), \qquad \bx\in \R^d_+.
$$
\end{lemma}

\begin{proof}
We slightly modify the argument in Resnick \cite{Resnick:2008}, p.~267-268, to account for the fact that $\alpha$ need not be one.

Define the polar coordinate transformation $T:\bx\mapsto (\|\bx\|_1,\bx/\|\bx\|_1)$. Then by Proposition \ref{prop:mv_regvar}, $\mu\circ T^{-1}(r,A)=r^{-\alpha}\sigma(A)$ for $r>0$, $A\subseteq\mathcal S^d$, with density $\alpha r^{-\alpha-1}dr\sigma(d\boldsymbol w)$. Of course $\mu([\0,\bx]^c) = \mu\circ T^{-1} \circ T([\0,\bx]^c)$, and the same computations as in \cite{Resnick:2008}, p.~267 yield
$$
T([\0,\bx]^c) = \left\{ (r,\boldsymbol w) \mid r > \min\left\{\frac{x_1}{w_1},\ldots,\frac{x_d}{w_d}\right\} \right\}.
$$
We can now integrate the density of $\mu\circ T^{-1}$ over this region:
\begin{eqnarray*}
\mu([\0,\bx]^c) &=& \int \int_{T([\0,\bx]^c)} \alpha r^{-\alpha-1} dr \sigma(d\boldsymbol w) \\
&=& \int_{\mathcal S^d}  \int_{ \left\{ r > \min\left\{\frac{x_1}{w_1},\ldots,\frac{x_d}{w_d}\right\} \right\} } \alpha r^{-\alpha-1} dr \sigma(d\boldsymbol w) \\
&=& \int_{\mathcal S^d} \left[ \min\left\{\frac{x_1}{w_1},\ldots,\frac{x_d}{w_d}\right\} \right]^{-\alpha} \sigma(d\boldsymbol w),
\end{eqnarray*}
which is equal to the asserted expression.
\proofend \end{proof}

It is now possible to compute the copula $C_0$ of the limiting distribution $G$. Although it does not have a particularly elegant form (except in the bivariate case, as we shall see), it does have an important interpretation: if $C$ is the copula of the distribution $\boldsymbol Y$, then corresponding $n$-block \emph{minima} have a dependence structure approximated by $C_0$ for large $n$. Equivalently, $C_0$ describes the dependence of the $n$-block \emph{maxima} of $\boldsymbol Y$ in case $C$ is the \emph{survival copula} of $\boldsymbol Y$. See the discussion at the end of Section \ref{section:L1nsd_min}.

First, to obtain the marginals $G_1$ of $G$ (they are all equal due to exchangeability) we note that by Proposition \ref{prop:mv_regvar}, $G_1(x)=\exp\{-\mu(B_x)\}$, where $B_x:=[\0,(x,\infty,\ldots,\infty)]^c$. Lemma \ref{lemma:mu_repr} then gives
$$
\mu(B_x)= x^{-\alpha} \int_{\mathcal S^d} w_1^\alpha \sigma(d\boldsymbol w) = c_\alpha x^{-\alpha},
$$
with $c_\alpha=\int_{\mathcal S^d} w_1^\alpha \sigma(d\boldsymbol w)$. Hence $G_1(x)=\exp\{-c_\alpha x^{-\alpha}\}$, which in particular means that $G_1$ is of the same type as $\Phi_\alpha$. This is to be expected since  the marginal survival function of $\bX$ is $\psi=\Wd F$, which is in $\MDA(\Phi_\alpha)$ by Proposition \ref{prop:Wd_tail_s1} and Theorem \ref{thm:MDA_char} (we still work under the assumption of Proposition \ref{prop:mv_regvar}: $R\in\RV_{-\alpha}$ with $\alpha\in(0,\infty)$).

Inverting $G_1$ gives $G_1^{\leftarrow}(u) = (-c_\alpha^{-1}\log u )^{-1/\alpha}$, and we may apply Lemma \ref{lemma:mu_repr} to compute the copula of $G$ as $C_0(\bu) = G(G_1^{-1}(u_1),\ldots,G_1^{-1}(u_d))$:
\begin{eqnarray*}
C_0(\boldsymbol u) &=& \exp\left\{ - \int_{\mathcal{S}^d} \max \left\{-\frac{\log u_1}{c_\alpha}w^\alpha_1,\ldots, -\frac{\log u_d}{c_\alpha}w^\alpha_d \right\}\sigma(d \boldsymbol{w}) \right\} \\
&=&  \exp\left\{ - \frac{1}{c_\alpha} \int_{\mathcal{S}^d} \max \left\{-w^\alpha_1 \log u_1,\ldots, -w^\alpha_d \log u_d  \right\}\sigma(d \boldsymbol{w}) \right\}.
\end{eqnarray*}

Although by Sklar's Theorem we already know that $C_0$ is a copula, one may readily verify that the marginals are standard uniform. It is also easy to check that $C_0$ satisfies the EV copula property $C_0^n(\boldsymbol u^{1/n})=C_0(\boldsymbol u)$.

\begin{example} \label{ex:gal}
For $d=2$, $C_0$ can be computed explicitly due to the simple form of the spectral measure $\sigma$. The first step is to calculate $V(x_1,x_2):=\int_{\mathcal S^2} \max\{w^\alpha_1 x_1, w^\alpha_2 x_2\} \sigma(d\boldsymbol w)$ for $x_1,x_2>0$. We parameterize $\mathcal S^2$ by $t:=w_1$ and obtain
$$
V(x_1,x_2) = \int_0^1 \max\{ x_1 t^\alpha, x_2 (1-t)^\alpha \} dt = \int_0^\kappa x_2 (1-t)^\alpha dt + \int_\kappa^1 x_1 t^\alpha dt,
$$
where $\kappa=x_2^{1/\alpha}/(x_1^{1/\alpha}+x_2^{1/\alpha})$ solves the equation $x_1 \kappa^\alpha = x_2 (1-\kappa)^\alpha$. The two integrals on the right side are straightforward to evaluate, and the result is
$$
V(x_1,x_2) = \frac{1}{1+\alpha} \left(x_1+x_2 - \left( x_1^{-1/\alpha} + x_2^{-1/\alpha}\right)^{-\alpha} \right).
$$
The quantity $c_\alpha$ is obtained similarly: $c_\alpha=\int_{\mathcal S^d} w^\alpha_1 \sigma(d\boldsymbol w) = \int_0^1 t^\alpha dt = (1+\alpha)^{-1}$. Piecing things together we get $C_0(u_1,u_2) = \exp\left\{ - (1+\alpha) V(-\log u_1, -\log u_2) \right\}$, or, more explicitly,
$$
C_0(u_1,u_2) = u_1 u_2 \exp\left\{ \left( (-\log u_1)^{-1/\alpha} + (-\log u_2)^{-1/\alpha} \right)^{-\alpha} \right\}.
$$
This is nothing but the well-known Galambos copula.
\end{example}

We conclude this section with a remark regarding the fact that we only considered the heavy-tailed case, i.e.~$R\in\RV_{-\alpha}$ for $\alpha\in(0,\infty)$, or, equivalently, $R$ being in the Fr\'echet domain of attraction. What would happen if $R$ were in the Gumbel domain? Assuming unbounded support for simplicity we have
$$
\lim_{u\rightarrow\infty} \frac{\P[\|\bX\|_1>ux, \bX/\|\bX\|_1 \in A]}{\P[\|\bX\|_1>u]} = \lim_{u\rightarrow\infty} \frac{\bar{F}(ux)}{\bar{F}(u)}\P[\boldsymbol{S}\in A] = 0
$$
for all Borel sets $A\subseteq\mathcal S^d$. The limiting measure $\mu$ in Proposition \ref{prop:mv_regvar} would thus be identically zero, $\mu([\0,\bx]^c)=0$ for $\bx\in \R^d_+\backslash \{\0\}$, meaning that the limit distribution would be degenerate with a point mass at the origin. Although a different normalization (such as $(\cdot-u)/a(u)$ as in Theorem \ref{thm:MDA_char}) might still yield meaningful results, we do not pursue this here since it involves translations and thus is non-standard in the basic theory of multivariate regular variation and multivariate extremes.

\chapter{Archimedean threshold copulas} \label{section:threshold}

In the univariate theory of extremes, the problem of determining the limiting distributions of maxima (or minima) is intimately related to that of finding the limiting distributions of threshold exceedances. In the multivariate setting things become more complicated, since thresholds can have many different forms and diverge in many different ways. In this chapter we consider the limiting dependence structures of \emph{componentwise} exceedances when the underlying distribution has an Archimedean copula. The  crucial tool will, as before, be the $\ell_1$-norm symmetric distributions.

\section{Lower threshold copulas}

Let $C$ be a $d$-dimensional copula and suppose $(U_1,\ldots,U_d)$ is distributed according to $C$. For given quantile levels $\alpha_i\in(0,1]$,
\begin{equation} \label{eqn:LLTC1}
\P[U_1 \leq \beta_1, \ldots, U_d\leq\beta_d \mid U_1 \leq \alpha_1,\ldots, U_d\leq \alpha_d], \qquad \beta_i\in[0,\alpha_i],
\end{equation}
defines a distribution function on $[0,\alpha_1]\times\ldots\times[0,\alpha_d]$. The copula of this distribution is called the \emph{lower threshold copula} of $C$; the limit as each $\alpha_i$ decreases is called the \emph{limiting lower threshold copula} (LLTC). If $C$ is Archimedean, the LLTC is always the Clayton copula, whose generator is $\psi^{\mathrm{Cl}}_\theta(x)=(1+\theta x)^{-1/\theta}$. This result follows in a very transparent way using the connection with $\ell_1$-norm symmetric distributions.

We need to be more precise about the statement ``each $\alpha_i$ decreases''. It can happen that the set $L_0:=\{\bu\in[0,1]^d \mid C(\bu)=0\}$, called the \emph{zero set} of $C$, has positive measure (w.r.t.~Lebesgue measure on $[0,1]^d$.) If $C$ is Archimedean, this happens when the generator $\psi$ is non-strict. In this case one cannot let $\alpha_i\downarrow0$, since eventually the probabilities in (\ref{eqn:LLTC1}) will not be well-defined. Instead, we should require $\boldsymbol{\alpha}:=(\alpha_1,\ldots,\alpha_d)\notin L_0$ and $\boldsymbol{\alpha}\rightarrow\bu_0$ for some $\bu_0\in \partial L_0$. If $L_0\neq\{\0\}$ there are different possible choices for $\bu_0$, and even when $L_0=\{\0\}$ there are no guarantees a priori that the LLTC is independent of how $\boldsymbol{\alpha}$ approaches $\bu_0=\0$. To avoid ambiguities, one therefore usually takes $\alpha_1=\alpha_2=\ldots=\alpha_d$ and lets $\alpha_i$ decrease until it hits $\partial L_0$ (this is the case in e.g.~McNeil et al. \cite{McNeil/Frey/Embrechts:2005} and Juri and W\"uthrich \cite{Juri/Wuethrich:2002}). Since in the Archimedean case, $L_0=\{\bu\in[0,1]^d \mid \sum_{i=1}^d \psi^{-1}(u_i)\geq\psi^{-1}(0) \}$, we are able to handle this complication by supposing that $\sum_{i=1}^d\psi^{-1}(\alpha_i) \uparrow \psi^{-1}(0)$, which ensures that $\boldsymbol{\alpha}$ approaches $L_0$ from above. Moreover, it turns out that neither the relative rates of convergence of the $\alpha_i$, nor which $\bu_0\in L_0$ they tend to, affects the LLTC. We therefore allow the quantiles $\alpha_i$ to be different in order to emphasize this fact.

Now, suppose $C$ is Archimedean with generator $\psi$. Since $\psi^{-1}$ is decreasing, $\{U_1\leq u\} = \{ \psi^{-1}(U_1)\geq \psi^{-1}(u)\}$ for $u\in [0,1]$. With $X_i:=\psi^{-1}(U_i)$, $v_i:=\psi^{-1}(\alpha_i)$ and $x_i:=\psi^{-1}(\beta_i)-\psi^{-1}(\alpha_i)\geq 0$, (\ref{eqn:LLTC1}) can thus be written as
\begin{equation} \label{eqn:Hbar1}
\left. \begin{array}{l}
\vspace{3mm}
\P[U_1 \leq \beta_1, \ldots, U_d\leq\beta_d \mid U_1 \leq \alpha_1,\ldots, U_d\leq \alpha_d] =  \\
\vspace{3mm}
\qquad \qquad = \P[X_1-v_1 \geq x_1, \ldots, X_d-v_d \geq x_d \mid X_1 \geq v_1, \ldots, X_1 \geq v_1] \\
\qquad \qquad =: \bar{H}_{\boldsymbol{v}}(\bx),
\end{array} \right.
\end{equation}
where $\bx=(x_1,\ldots,x_d)$ with $x_i\in[0,\psi^{-1}(0)-v_i]$ and $\boldsymbol{v}=(v_1,\ldots,v_d)$ denotes the vector of thresholds. Since all $X_i$ are continuous, the inequalities can be taken as strict, and we see that $\bar H_{\boldsymbol v}$ defines a survival function on the interior of $\R^d_+$.

\begin{lemma} \label{lemma:Hbar_L1-nsd}
$\bX=(X_1,\ldots,X_d)$ as defined above is $\ell_1$-norm symmetric with survival copula generator $\psi$. Moreover, $\bX_{\boldsymbol{v}}\sim \bar{H}_{\boldsymbol{v}}$ is also $\ell_1$-norm symmetric with survival copula generator $\psi_{\boldsymbol v}$ given by
$$
\psi_{\boldsymbol{v}}(x) = \frac{\psi(x+\|\boldsymbol{v}\|_1)}{\psi(\|\boldsymbol{v}\|_1)}.
$$
\end{lemma}

\begin{proof}
The statement about $\bX$ follows from Corollary \ref{cor:main}. An explicit computation yields
$$
\bar{H}_{\boldsymbol{v}}(\bx) = \frac{\P[X_1-v_1 > x_1, \ldots, X_d-v_d > x_d]}{\P[X_1 > v_1, \ldots, X_1 > v_1]} = \frac{\psi(\|\bx\|_1+\|\boldsymbol{v}\|_1)}{\psi(\|\boldsymbol{v}\|_1)},
$$
whence the last assertion due to Theorem \ref{thm:main}.
\proofend \end{proof}

Due to (\ref{eqn:Hbar1}), the survival copula corresponding to $\bar{H}_{\boldsymbol{v}}$, which we denote by $C_{\boldsymbol{v}}$, is the lower threshold copula of $C$. Now, $\alpha_i \downarrow 0$ if and only if $v_i\uparrow \psi^{-1}(0)$, so in view of the discussion in the beginning of this section the LLTC of C is given by the limit of $C_{\boldsymbol{v}}$ as $\|\boldsymbol{v}\|_1\uparrow\psi^{-1}(0)$.

\begin{proposition} \label{prop:L1-nsd_thr_conv}
Let $\psi$ and $\bX_{\boldsymbol v}$ be as in Lemma \ref{lemma:Hbar_L1-nsd}, and suppose $\psi$ is in some domain of attraction. Then there exists a normalization $\beta(\cdot)$ such that $\boldsymbol X_{\boldsymbol v}/\beta(\|\boldsymbol v\|_1)$ converges weakly to an $\ell_1$-norm symmetric distribution with survival copula generator $\psi^{\mathrm{Cl}}_\theta$ as $\|\boldsymbol{v}\|_1\rightarrow\psi^{-1}(0)$. The parameter $\theta$ is given by:
\begin{itemize}
\item[(i)] If $\psi$ is in the Fr\'echet domain, i.e.~$\psi(x)\in \RV_{-\alpha}$, for some $\alpha\in(0,\infty)$, then $\theta=1/\alpha$.
\item[(ii)] If $\psi$ is in the Gumbel domain, then $\theta=0$.
\item[(iii)] If $\psi$ is in the Weibull domain, i.e.~$\psi(x^*-x^{-1})\in \RV_{-\alpha}$ for some $\alpha\in(0,\infty)$ and where $x^*=\psi^{-1}(0)$, then $\theta=-1/\alpha$.
\end{itemize}
\end{proposition}

\begin{proof}
In the Fr\'echet case Lemma \ref{lemma:Hbar_L1-nsd} and the fact that $\psi\in\RV_{-\alpha}$ yield
$$
\P\left[\frac{\boldsymbol X_{\boldsymbol v}}{\|\boldsymbol v\|_1} > \boldsymbol x\right] = \frac{\psi(\|\boldsymbol{v}\|_1(1 + \|\bx\|_1))}{\psi(\|\boldsymbol{v}\|_1)} \to (1+\|\boldsymbol x\|_1)^{-\alpha} = \psi^{\mathrm{Cl}}_{1/\alpha}(\|\boldsymbol x\|_1)
$$
as $\|\boldsymbol v\|_1\to \psi^{-1}(0)=\infty$. In this case we may thus choose $\beta(x)=x$. In the Weibull case, let $t$ be defined via $\|\boldsymbol v\|_1 = x^*-t^{-1}$, and set $\beta(x):=x^*-x$. Then
\begin{eqnarray*}
\P\left[\frac{\boldsymbol X_{\boldsymbol v}}{\beta(\|\boldsymbol v\|_1)} > \boldsymbol x\right] &=& \frac{\psi\left(\|\boldsymbol{v}\|_1 + \beta(\|\boldsymbol v\|_1)\|\bx\|_1\right)}{\psi(\|\boldsymbol{v}\|_1)} \\
&=& \frac{\psi(x^* - (1-\|\boldsymbol x\|_1)t^{-1})}{\psi(x^*-t^{-1})} \\
&\to& (1-\|\boldsymbol x\|_1)^{\alpha} = \psi^{\mathrm{Cl}}_{-1/\alpha}(\|\boldsymbol x\|_1), \qquad \|\boldsymbol v\|_1\to x^*,
\end{eqnarray*}
since $\psi(x^*-x^{-1})\in\RV_{-\alpha}$. In the Gumbel case we take $\beta$ to be the auxiliary function of $\psi$. The right side of the first line in the previous display now tends to $\exp\{-\|\boldsymbol x\|_1\} = \psi^{\mathrm{Cl}}_0(\|\boldsymbol x\|_1)$ due to the property (\ref{eqn:gamma}) in Theorem \ref{thm:MDA_char}. In each of the three cases, the assertion that the limit is $\ell_1$-norm symmetric follows from Theorem \ref{thm:main}.
\proofend \end{proof}

\textbf{Remark.}
The function $\beta(\cdot)$ has the standard interpretation as the proper normalization for the conditional univariate exceedances $(X_1-v|X_1>v)$ to converge to a Generalized Pareto Distribution (cf. Balkema and de Haan \cite{Balkema/deHaan:1974} and Pickands \cite{Pickands:1975}).

Note that the parameter $\theta$ can be negative. This happens when $a:=\psi^{-1}(0)<\infty$ and $1-\psi(a-x^{-1})\in RV_{-\alpha}$, $\alpha>0$. From Proposition \ref{prop:nonstrict_regvar} we know that one has the limitation $\theta=-1/\alpha\geq-1/(d-1)$, which is consistent with the range of parameter values for which the Clayton copula is defined (see McNeil and Ne\v{s}lehov\'a \cite{McNeil/Neslehova:2008}, Example 4.)

\begin{corollary} \label{cor:LLTC_limit}
Let $C$ be an Archimedean copula whose generator $\psi$ is in some domain of attraction. The corresponding limiting lower threshold copula (LLTC) is then Clayton with parameter $\theta$ as in Proposition \ref{prop:L1-nsd_thr_conv}.
\end{corollary}

\begin{proof}
Denoting again the copula if $\bar H_{\boldsymbol v}$ by $C_{\boldsymbol v}$, the LLTC is given by the pointwise limit of $C_{\boldsymbol{v}}$ as $\|\boldsymbol{v}\|_1\rightarrow\psi^{-1}(0)$. But scaling by a positive constant does not change the copula\footnote{In fact, copulas are invariant under strictly increasing transformations $T_i$ of the components of the underlying random vector---if $(Y_1,Y_2)\sim H$ with marginals $H_1$, $H_2$ has copula $C$, then the copula of $(T_1(Y_1),T_2(Y_2))$ is $H(T_1^{\leftarrow}\circ T_1 \circ H_1^{\leftarrow}(\cdot), T_2^{\leftarrow}\circ T_2 \circ H_2^{\leftarrow}(\cdot))=C$ by the properties of generalized inverses.}, so $\bX_{\boldsymbol{v}}/\beta(\|\boldsymbol{v}\|_1)$, like $\bX_{\boldsymbol{v}}$, has survival copula $C_{\boldsymbol{v}}$. Since weak convergence of the distributions implies pointwise convergence of the copulas (cf.~Feidt et al.~\cite{Feidt/Genest/Neslehova:2008}), the assertion follows from Proposition \ref{prop:L1-nsd_thr_conv}. 
\proofend \end{proof}

The fact that only the sum $\|\boldsymbol{v}\|_1$ has to tend to $\psi^{-1}(0)$ immediately gives the form of all \emph{marginal LLTCs}, which are defined as the limits of
\begin{equation} \label{eqn:marg_LLTC}
\P[U_i \leq \beta_i, \ i\in J \mid U_i \leq \alpha_i, \ i\in J], \qquad \beta_i\in[0,\alpha_i], 
\end{equation}
as $(\alpha_i)_{i\in J}$ tends to the corresponding zero set, where $J\subset\{1,\ldots,d\}$ is a non-empty index set. More precisely, we have

\begin{corollary}
Let $C$ be an Archimedean copula whose generator $\psi$ is in some domain of attraction. All corresponding marginal LLTCs are then Clayton with parameter $\theta$ as in Proposition \ref{prop:L1-nsd_thr_conv}.
\end{corollary}

\begin{proof}
(\ref{eqn:marg_LLTC}) is obtained from (\ref{eqn:LLTC1}) upon setting $\alpha_i=1$ and $\beta_i=1$ for all $i\notin J$. This corresponds to $x_i=v_i=0$, $i\notin J$. Letting the remaining $v_i$ increase such that $\|\boldsymbol{v}\|_1\uparrow\psi^{-1}(0)$, the assertion is inferred from Proposition \ref{prop:L1-nsd_thr_conv} as in Corollary \ref{cor:LLTC_limit}. 
\proofend \end{proof}

In a modeling situation, there might be given requirements on the behavior of the LLTC of the copula $C$. If one is working with the stochastic representation in terms of $\ell_1$-norm symmetric distributions, the results in Chapter \ref{section:Wd_tail} translate the three conditions on $\psi$ given in Proposition \ref{prop:L1-nsd_thr_conv} into conditions on the corresponding radial part $R$. It may also be the case that a certain $R$ has been chosen through some inference procedure, and one then wants to know the implications for the extremal behavior of the limiting copula. Again, Chapter  \ref{section:Wd_tail} provides the necessary link.

We conclude with a remark on the connection between threshold exceedances and the componentwise maxima treated in Section \ref{section:L1nsd_max}. One may view the conclusion of Example \ref{ex:gal} as an extension of the well-known fact that the Galambos copula is the EV limit of the survival copula of the bivariate Pareto distribution (see McNeil et al.~\cite{McNeil/Frey/Embrechts:2005}, Example 7.50). In fact, the (symmetric) bivariate Pareto is $\ell_1$-norm symmetric and its survival copula is Clayton (this is easily checked: its joint survival function is $\bar H(x_1,x_2)=(x_1+x_2+1)^{-\alpha}$ with $\alpha>0$).

Now, by Proposition \ref{prop:L1-nsd_thr_conv}, any heavy-tailed bivariate $\ell_1$-norm symmetric distribution has asymptotically the same shape as the bivariate Pareto far out in the tail. In the case of $\ell_1$-norm symmetric distributions, therefore, the following recipe produces correct results: given $\bX=(X_1,X_2)$, compute the limiting distribution $\bX_\infty$ of the threshold exceedances. Take the componentwise maxima of $\bX_\infty$ and compute the corresponding limiting distribution. The copula of this limit is the same as the copula of the distribution that is obtained by directly taking the limit of the componentwise maxima of $\bX$. Intuitively such a result is reasonable; for large $n$, the $n$-block maxima will be located far out in the tail, where the underlying distribution looks approximately like the limit of the exceedances.

The last statement gives rise to an interesting question: in general, what are the conditions on the underlying multivariate distribution for these heuristics to be correct? I.e., under which circumstances does the recipe described above give correct results?

\section{Upper threshold copulas}

In analogy with the lower threshold copulas, one may define the \emph{upper threshold copula} for given quantiles $\alpha_1,\ldots,\alpha_d$ as the copula of the distribution
\begin{equation} \label{eqn:LUTC1}
\P[U_1 > \beta_1, \ldots, U_d>\beta_d \mid U_1 > \alpha_1,\ldots, U_d> \alpha_d], \qquad \beta_i\in[\alpha_i,1],
\end{equation}
where again $(U_1,\ldots,U_d)$ are distributed according to the copula $C$. The limit as the $\alpha_i\uparrow1$ is called the \emph{limiting upper threshold copula} (LUTC)\footnote{We do not have the same problem as in the LLTC case---the set $L_1:=\{\bu\in[0,1]^d\ |\ C(\bu)=1\}$ is always a nullset. Indeed, $C(u_1,\ldots,u_d)\leq C(u_1,1,\ldots,1)=u_1$ implies that $C(\boldsymbol u)<1$ for $\boldsymbol u \neq \boldsymbol 1$.}, and our aim is to investigate the extent to which our previous methodology can be applied to compute the LUTC in the case when $C$ is Archimedean.

We first define $X_i:=\psi^{-1}(U_i)$, $v_i:=\psi^{-1}(\alpha_i)$ and $x_i:=\psi^{-1}(\beta_i)$. The quantity (\ref{eqn:LUTC1}) can then be expressed in terms of the $\ell_1$-norm symmetric random variable $\bX=(X_1,\ldots,X_d)$, whose survival copula is $C$ (that $\bX$ is $\ell_1$-norm symmetric follows as before from Corollary \ref{cor:main}):
\begin{equation*}
\left. \begin{array}{l}
\vspace{3mm}
\P[U_1 > \beta_1, \ldots, U_d>\beta_d \mid U_1 > \alpha_1,\ldots, U_d> \alpha_d] = \\
\vspace{3mm}
\qquad \qquad = \P[\bX\leq\bx \mid \bX \leq \boldsymbol{v}] \\
\qquad \qquad =: H_{\boldsymbol{v}}(\bx), \qquad \qquad x_i\in[0,v_i].
\end{array} \right.
\end{equation*}
To get the inequalities right we used strict monotonicity of $\psi^{-1}$ and the fact that $\bX$ has continuous marginals. It is straightforward to show that $H_{\boldsymbol{v}}$ is \emph{not} the distribution function of a $\ell_1$-norm symmetric distribution---indeed, the support of $(\bX \mid \bX \leq \boldsymbol{v})$ is $[\0,\boldsymbol v]$, which is not consistent with $\ell_1$-norm symmetry. Its survival copula is therefore not Archimedean, despite the fact that $C$ is. Hence the methodology used in the LLTC case cannot be employed, and this difficulty cannot be overcome using basic multivariate extreme value theory as easily as in the EV copula case in Section \ref{section:L1nsd_max}. In fact, the ULTC cannot be explicitly computed, but we shall derive a distribution function whose survival copula is the ULTC of $C$. For simplicity we carry out the calculations in two dimensions to illustrate the procedure. Generalizing the results to higher dimension is straightforward.

With $d=2$, we have
\begin{eqnarray*}
\P[X_1\leq x_1,X_2\leq x_2] &=& 1- \P[X_1>x_1] - \P[X_1>x_1] + \P[X_1> x_1,X_2> x_2] =\\
&=& 1-\psi(x_1)-\psi(x_2)+\psi(x_1+x_2) = \\
&=& G(x_1)+G(x_2)-G(x_1+x_2),
\end{eqnarray*}
where $G(x):=1-\psi(x)$ is the marginal distribution function of $\bX$. In the general case one may use the sieve formula to express $\P[\bX\leq \bx]$ in terms of $G(\sum_{i\in J} x_i)$ for subsets $J\subseteq\{1,\ldots,d\}$. Proceeding with the case $d=2$ we obtain
$$
H_{\boldsymbol{v}}(\bx) = \frac{G(x_1)+G(x_2)-G(x_1+x_2)}{G(v_1)+G(v_2)-G(v_1+v_2)}, \qquad x_i\in[0,v_i].
$$
Now, $\alpha_i\uparrow1$ is equivalent to $v_i\downarrow0$, so in order to have a non-degenerate limit we define $\tilde{\bX}=(\tilde{X}_1,\tilde{X}_2)$ as the normalized $[0,1]^2$-valued random variable (depending on $\boldsymbol{v}$) distributed according to
\[ \P[\tilde{X}_1\leq x_1,\tilde{X}_2\leq x_2] = H_{\boldsymbol{v}}(x_1v_1, x_2v_2), \qquad x_i\in[0,1]. \]
The copula is not affected, since it is invariant under strictly increasing transformations of the components. We may therefore consider the copula of $\tilde{\bX}$ as $v_1$ and $v_2$ tend to zero.

\begin{proposition} \label{prop:7.5}
Suppose $1/X_1 \sim \psi(x^{-1})$ is in some domain of attraction, and let $v_1=v_2=v$. In particular we then have $1-\psi(x^{-1})\in RV_{-\alpha}$ for some $\alpha\in(0,1]$, and if in addition $\alpha<1$, then
\[ H_0(x_1,x_2) := \lim_{v\downarrow0} \P[\tilde{X}_1\leq x_1,\tilde{X}_2\leq x_2]  = \frac{x_1^\alpha+x_2^\alpha - (x_1+x_2)^\alpha}{2-2^\alpha}. \]
\end{proposition}

\begin{proof}
Theorem \ref{thm:MDA_char} implies that $G(x^{-1})=1-\psi(x^{-1})\in RV_{\alpha}(0)$ for some $\alpha\in(0,\infty]$, and by Proposition \ref{prop:2}, $\alpha\leq 1$. Since $v_1=v_2=v$ this yields
\begin{eqnarray*}
\P[\tilde{X}_1\leq x_1,\tilde{X}_2\leq x_2] &=& \frac{G(x_1v)+G(x_2v)-G(x_1v+x_2v)}{2 G(v)-G(2v)} = \\
&=& \frac{ \frac{G(x_1v)}{G(v)} +  \frac{G(x_2v)}{G(v)} - \frac{G(v(x_1+x_2))}{G(v)}  }{ 2 - \frac{G(2v)}{G(v)}} \\
&\rightarrow&\frac{x_1^\alpha+x_2^\alpha - (x_1+x_2)^\alpha}{2-2^\alpha}, \qquad v\downarrow0.
\end{eqnarray*}
\proofend \end{proof}

The assumption $v_1=v_2$ is not restrictive if we are only interested in the copulas. Indeed, if $v_2/v_1\rightarrow c>0$, then one can show, using that the convergence in Proposition \ref{prop:7.5} is locally uniform (a basic result from Karamata theory), that $H_0(x_1,cx_2)$ is obtained in the limit. This just represents a linear scaling of the second component, which leaves the copula unchanged. The marginals $H_0(x,1)=(x^\alpha+1-(x+1)^\alpha)/(2-2^\alpha)$ do not have closed form inverses, so we are not able to explicitly compute the copula or the survival copula. As a final remark we mention that more precise results can be obtained in the limiting case when $\alpha=1$, see Charpentier and Segers \cite{Charpentier/Segers:2008}.

\bibliographystyle{plain}

%\bibliography{LN}

\end{document}